\documentclass{amsart}
\usepackage{amsmath}
\usepackage{amscd,amssymb}
\usepackage{amsmath,amsfonts,amsthm}

\usepackage{bbm}

\newcommand{\R}{\mathbb{R}}

\newcommand{\ob}[2]{{\bf U}({#1},{#2})}

\newcommand{\altdown}[1]{\bigwedge\nolimits_{#1}}
\newcommand{\op}[2]{{\bf #1}\,#2}

\newcommand{\ul}[1]{\underline{#1}}

\newcommand{\Lip}[1]{{\rm Lip}\left(#1\right)}

\newcommand{\Tan}{{\bf Tan}}
\newcommand{\Nor}{{\bf Nor}}

\newcommand{\Lin}{\mathbbm{Lin}}
\newcommand{\End}{\mathbbm{End}}

\newcommand{\MultiLin}{\mathbbm{MultiLin}}

\newcommand{\N}{\mathbb{N}}

\newcommand{\cA}{\mathcal{A}}
\newcommand{\cB}{\mathcal{B}}
\newcommand{\cC}{\mathcal{C}}

\newcommand{\cE}{\mathcal{E}}

\newcommand{\cG}{\mathcal{G}}
\newcommand{\cH}{\mathcal{H}}

\newcommand{\cJ}{\mathcal{J}}

\newcommand{\cL}{\mathcal{L}}
\newcommand{\cM}{\mathcal{M}}
\newcommand{\cN}{\mathcal{N}}

\newcommand{\cP}{\mathcal{P}}

\newcommand{\cS}{\mathcal{S}}

\newcommand{\cU}{\mathcal{U}}
\newcommand{\cV}{\mathcal{V}}
\newcommand{\cW}{\mathcal{W}}

\newcommand{\bB}{{\bf B}}

\newcommand{\bE}{{\bf E}}

\newcommand{\bJ}{{\bf J}}

\newcommand{\bL}{{\bf L}}

\newcommand{\bS}{{\bf S}}

\newcommand{\bc}{{\bf c}}

\newcommand{\bj}{{\bf j}}

\newcommand{\bo}{{\bf o}}

\newcommand{\br}{{\bf r}}

\newcommand{\bz}{{\bf z}}

\def\altdown#1{\bigwedge\nolimits_{#1}}

\usepackage{amsthm}
\theoremstyle{plain}
\newtheorem{thm}{Theorem}[section]
\newtheorem{prop}{Proposition}[section]
\newtheorem{lem}{Lemma}[section]
\newtheorem{cor}{Corollary}[section]

\newtheorem{rem}{Remark}[section]
\newtheorem*{rem*}{Remark}
\newtheorem{defn}{Definition}[section]
\newtheorem*{defn*}{Definition}
\newtheorem{ex}{Example}[section]

\newcommand{\Sm}{\mathbb{S}^{n-1}}

\newcommand{\dist}{{\bf dist}}

\newcommand{\bbeta}{\boldsymbol\eta}

\newcommand{\ubr}{\underline{\bold r}}

\newcommand{\ubb}{\underline{\bold b}}

\newcommand{\reach}{{\bf reach}}

\begin{document}

\title
{Corrections to a paper of Allard and Almgren on the uniqueness of tangent cones}

\author{William K. Allard}
\thanks{Supported in part by Los Alamos National Laboratory.}

\date{\today}

\maketitle

\setcounter{tocdepth}{3}
\tableofcontents

\newcommand{\mO}{\mathbb{O}}
\newcommand{\mG}{\mathbb{G}}
\newcommand{\mGO}{\mathbb{GO}}

\newcommand{\Sob}[3]{{\bf U}_{#1}(#2,#3)}\newcommand{\Scb}[3]{{\bf B}_{#1}(#2,#3)}



\section{Introduction.}

The paper referred to in the title is \cite{AA81}.

Several months ago Francesco Maggi emailed me saying that the inequalities
 5.3(4),(5) of \cite{AA81} were wrong. In fact, as he pointed out,
 their incorrectness is immediately apparent if one takes $Z=0$ there.
 Maggi and his coauthor wanted to use these inequalities in their paper
 \cite{MN}. They were able to obtain a version of these inequalities
 which suffice for the carrying out the work in \cite{MN}.

 I started writing this paper in order to provide a 
 version of these inequalities as needed in [AA81].
In thinking about this material I began to realize there were other problems
with the paper. As a result I ended up {\em completely rewriting 5.1-5.4 on pages 
243-248 of \cite{AA81}}; this rewrite is the contents of this paper. In addition to many annoying misprints many of the
needed definitions and proofs in 5.1-5.4 are incomplete or absent. For example,
it is not said where $z$ in 5.1(2) comes from; this omission completely
surprised me since I remember doing a lot of work to come up with $z$
when \cite{AA81} was being written. Also, much of the necessary material about
the reach of a submanifold as in \cite{FE2} is not provided in \cite{AA81}.

The table of contents can serve as an index. In particular one sees there
where the constants $\epsilon_1$ through $\epsilon_6$ are introduced.

This material is extremely technical. One way to navigate this paper would
be to start looking at Proposition \ref{Box} and work backwards.

\section{Preliminaries.}

If $X$ is a set we let
$$i_X=\{(x,x):x\in X\};$$
thus $i_X:X\rightarrow X$ and $i_X(x)=x$ for $x\in X$; we call $i_X$ the
{\bf identity map of $X$}; if $X\subset Y$ the $i_X$ is also called the
{\bf inclusion map of $X$ into $Y$}.

If $(X,\rho)$ and $(Y,\sigma)$ are metric spaces and $f:X\rightarrow Y$
we let
$$\op{Lip}{f}=\sup\left\{{\sigma(f(x_1),f(x_2))\over\rho(x,y)}:x_i\in X, \ i=1,2, \ \hbox{and} \ x_1\neq x_2\right\}.$$

\subsection{Some linear algebra.}

be the vector space of linear maps from $V$ to $W$. We let
$$\End(V)=\Lin(V,V).$$
If $V_1,\ldots,V_N$ and $W$ are vector spaces we let
$$\MultiLin(V_1,\ldots,V_N,W)$$
be the vector space of multilinear maps from $V_1\times\cdots\times V_N$
into $W$.

If $V$ and $W$ are vector spaces we let
$$\Lin(V,W)=\MultiLin(V,W).$$

If $V$ is a vector space we let
$$\End(V)=\Lin(V,V).$$

A {\bf Euclidean space} is a finite dimensional inner product space.
If $E$ is a Euclidean space and $l\in\N$ we let
$$\mG_l(E) \  \hbox{and} \ \mO_l(E)$$
be, respectively, the family of $l$ dimensional linear subspaces of $E$ and the family
of orthogonal projections of $E$ onto a member of $\mG_l(E)$.

If $E_i$, $i=1,\ldots,m$ and $F$ are Euclidean spaces and
$\mu:E_1\times\cdots\times E_m\rightarrow F$ is multilinear we let
$$||\mu||=\sup\{|\mu(v_1,\ldots,v_m)|:v_i\in E_i \ \hbox{and} \ |v_i|\leq 1,
 \ i=1,\dots,m\};$$
the vector space of multilinear maps from $E_1\times E_m$ into $F$ has a natural Euclidean structure and the norm with respect to this structure of $\mu$
as above will be larger than $||\cdot||$ unless the $E_i$ are one dimensional.
For example, if $i$ is the identity map of the Euclidean space $E$ then
$||i||=1$ but the Euclidean norm of $i$ is the square root of the dimension of
$E$.

\begin{prop}
Suppose  $S\in\mG_m(\R^n)$, $L\in\Lin(S,\R^n)$, $\op{ker}{L}=\{0\}$.
and $T=\op{rng}{L}$. and $\op{dim}{T}=m$.

There are unique $A$ and $O$ such that
$A\in\End(T)$, $O\in\Lin(S,T)$, $A^*=A$, $A(v)\bullet v> 0$ for $v\in T$, $A^2=L\circ L^*$, $O$ is an isometry and $L=A\circ O$.
\end{prop}
\begin{proof} $L\circ L^*$ is a positive definite and carries $T$
isomorphically onto $T$. Let $A$ be the unique positive definite
square root of $L\circ L^*$ and let $O=A^{-1}\circ L$.
\end{proof}

Suppose $S$ and $T$ are finite dimensional inner product spaces.
For $L,M\in\Lin(S,T)$
$$L\bullet M=\op{trace}{L^*\circ M}.$$
This is an inner product on $\Lin(S,T))$. If $u_1,\ldots,u_m$
is an orthonormal basis for $S$ we have
$$L\bullet M=\sum_{i=1}^m L(u_i)\bullet M(u_i).$$
It follows that
$$|L\bullet M|\leq m\,||L||\,||M||;$$
note that if $L$ is an isometry then $||L||=1$ and $|L|=\sqrt{m}$.

\subsection{A useful notational convention.}
Suppose $T$, $X$ and $Y$ are nonempty sets and $W\subset T\times X$.

If $h:T\times X\rightarrow Y$ and $t\in\cJ$ we let
$$h_t=\{(x,h(t,x)):(t,x)\in W\}.$$
Thus
$$h_t:\{x:(t,x)\in W\}\rightarrow Y\quad\hbox{for $t\in\cJ$};$$
of course $h_t=\emptyset$ if $W\cap(\{t\}\times X)=\emptyset$.

\section{Submanifolds of Euclidean spaces.}\label{subman}

\begin{defn}\label{mani} Suppose $E$ is a Euclidean space,
$m\in\N^+$ and $m\leq\op{dim}{E}$.
We let
$$\cM_m(E)$$
be the family of nonempty subsets $M$ of $E$ such that for each $a\in M$
there are $U,V,\psi,\phi$ such that $V$ is an open subset of $E$; $U$ is an open subset of $\R^n$,
$\psi:V\rightarrow U$, $\psi$ is smooth, $\phi:U\rightarrow V$, $\phi$
is smooth, $a\in V$ and
$$\psi(\phi(u))=u\quad\hbox{for $u\in U$.}$$
\end{defn}

\begin{rem} This a slight variation of the definition in [FE, 3.1.19(4).
  Note that $\phi=(\psi|M)^{-1}$
\end{rem}

A member of $\cM_m(E)$ is called a
{\bf smooth $m$ dimensional submanifold of $E$}. 

Evidently, $M\in\cM_{\op{dim}{E}}(E)$ if and only if $M$ is an open subset
of $E$. Also, if $L\subset M$ then $L\in\cM_m(E)$ if and only if $L$
is nonempty and is open relative to $M$.

It is elementary that
$$\Tan{M,\phi(u)}=\op{rng}{D\phi(u)}\in\mG_m(E)
\quad\hbox{for $u\in U$.}$$
It follows that $m$ is determined by $M$; we call it the {\bf dimension} of $M$.
We let
$$\tau:M\rightarrow\mG_m(E)
\quad\hbox{and}\quad
\nu:M\rightarrow\mG_{\op{dim}{E}-m}(E)$$
be the {\bf tangent map} and the {\bf normal map} of $M$, respectively;
by definition,
$\tau(x)$ is orthogonal projection of $E$ onto $\Tan{M,x}$ and $\nu(x)=\tau(x)^\perp$ for $x\in M$.

{\em We now fix $m,n\in\N^+$ with $m\leq n$.}

\subsection{}

For the remainder of this Section we fix $M\in\cM_m(\R^n)$
and a normed vector space $Y$.

\subsection{Smooth functions on $M$.}
A function $f:M\rightarrow Y$
is {\bf smooth} if $f\circ\phi$ is smooth if whenever $a,U,V,\phi,\psi$ are as in
Definition \ref{mani}.
$$\cE(M,Y)$$
be the vector space of smooth maps from $M$ into $Y$.

\begin{prop} Suppose $Y$ is a Euclidean space, $f:M\rightarrow Y$,
 $f$ is smooth and $a\in M$.  There is and only one
$$Df:M\rightarrow\Lin(\R^n,Y)$$
such that if $a\in M$ then
$$Df(a)(D\phi(\psi(a))(v))=D(f\circ\phi)(\psi(a))(v)
\quad\hbox{for $v\in\R^m$}$$
and
$$Df(a)\circ\nu(a)=0.$$

Moreover, if $N$ is a smooth submanifold of some Euclidean space
and $F:N\rightarrow M$ is smooth then $f\circ F$ is smooth and
$$D(f\circ F)(b)=Df(F(b))\circ Df(b)
\quad\hbox{for $b\in N$.}$$
\end{prop}
\begin{proof} Suppose $a,U,V,\phi,\psi$ are as in
Definition \ref{mani}. Then
$(f\circ F)|F^{-1}[V]=(f\circ\phi)\circ(\psi\circ F)|F^{-1}[V]$.
\end{proof}

If $f\in\cE(M,Y)$, $a\in M$ and $v\in E$ we let
$$D_v f(a)=Df(a)(v).$$
If $k\in\N^+$ and $v=(v_1,\ldots,v_k)\in E^k$ we let
$$D_v f(a)=D_{v_1}\cdots D_{v_k}f(a).$$

\begin{prop} Suppose $f\in\cE(M,Y)$, $a\in M$, $W$ is an open neighborhood
of $a$ in $\R^n$, $F:W\rightarrow Y$ and $F|M=f|W$. Then
$Df(a)(v)=DF(a)(v)$ whenever $a\in M\cap W$ and $v\in\Tan{M,a}$.
\end{prop}
\begin{proof} If $U,V,\phi,\psi$ are as in Definition \ref{mani} then $F\circ\phi$ equals $f\circ\phi$
in an open neighborhood of $\psi(a)$.
\end{proof}

We leave the straightforward proof of the following Proposition to the reader.
\begin{prop}\label{mani2}  Suppose $E$ is a Euclidean space,
$m\in\N^+$ and $m\leq\op{dim}{E}$, $M$ is a nonempty subset of $E$
and for each $a\in M$ there are $U,V,\psi,\phi$ such that $V$ is an open subset of $E$; $U$ is an open subset of some $m$ dimensional submanifold of some
Euclidean space,
$\psi:V\rightarrow U$, $\psi$ is smooth, $\phi:U\rightarrow V$, $\phi$
is smooth, $a\in V$ and
$$\psi(\phi(u))=u\quad\hbox{for $u\in U$.}$$

Then $M\in\cM_m(E)$.
\end{prop}

If $i_M$ is the inclusion map of $M$ in $\R^n$ then $i_M$ is smooth
and 
$$Di_M(x)=\tau(x)\quad\hbox{for $x\in M$.}$$
It follows that $\tau$ is smooth and since $\nu(x)=i_{\R^n}-\tau(x)$
for $x\in M$ we find that $\nu$ is smooth.

\begin{defn}\label{Df} For each $k\in\N$ we define
$$D^k f(a)\in\MultiLin(\underbrace{E,\ldots,E}_{\hbox{$k$ times}},Y)$$
by letting $D^0 f=f$, $D^1 f=Df$ and, if $k\geq 2$, requiring that
$$D^k f(a)(D\phi(\psi(a))(v_1),\ldots,D\phi(\psi(a))(v_k))
=D^k (f\circ\phi)(\psi(a))(v_1,\ldots,v_k)$$
whenever $a,U,V,\phi,\psi$ are as in Definition \ref{mani}
and $v_1,\ldots,v_k\in \R^m$
and that
$$D^k f(a)(w_1,\ldots,w_k)=0
\quad\hbox{whenever $w_1,\ldots,w_k\in E$ and
$w_i\in\Nor{M,a}$ for some $i=1,\ldots,m$.}$$
\end{defn}
As is well known $D^k f(a)$ is symmetric. We have
$$D^k f(a)(w_1,\ldots,w_k)=D_{w_1}\cdots D_{w_k} f(a)
\quad\hbox{whenever $w_1,\ldots,w_k\in E$.}$$

\begin{defn} Suppose $Y$ is a Euclidean space and
$f:M\rightarrow Y$ is smooth. For each $k\in\N$ we let
$$\chi^k(f)=\sup\{||D^k f(x)||:x\in M\}$$
and we let
$$\chi^{[k]}(f)=\max\{\chi^j(f):j\in \N \ \hbox{and} \ j\leq k\}.$$
Any of these numbers may be $\infty$.
\end{defn}

\begin{rem} We will {\em not} use $\chi$ as on page 229 of \cite{AA81}.
\end{rem}
 
\subsection{More on differentiation.}
Suppose $N\in\N^+$ and, for each $i=1,\ldots,N$, $E_i$ is a Euclidean space.
For each $i\in\{1,\ldots,N\}$ we let
$$\iota_i:E_i\rightarrow E_1\times\cdots\times E_N$$
be such that for each $u\in E_i$ and $j\in\{1,\ldots,N\}$
$$(\iota_i)_j(u)=\begin{cases}
  u&\hbox{if $j=i$,}\cr
  0&\hbox{if $j\neq i$.}
\end{cases}$$
For each $j\in\N^+$ we let
$$\iota_i^:E_i^j\rightarrow(E_1\times\cdots\times E_N)^j$$
be such that
$$\iota_i^j(u_1,\ldots,u_j)=(\iota_i(u_1),\ldots,\iota_i(u_j))
\quad\hbox{for $(u_1,\ldots,u_j)\in E_i^j$.}$$

Suppose $Y$ is a normed vector space and
$f:M_1\times\cdots\times M_N\rightarrow Y$ is smooth 
and $(k_1,\ldots,k_N)\in\N^N$ and $K>0$ where $K$ is the number
of positive $k_i$, $i=1,\ldots,N$.

Let $l_j$, $j=1,\ldots,K$, be such that
$$l_1<\cdots<l_K\quad\hbox{and}\quad\{l_j:j=1,\ldots,K\}=\{i\in\{1,\ldots,N\}:k_i>0\}.$$
We let
$$D^{k_1,\ldots,k_N}f:M_1\times\cdots\times M_N
\rightarrow\MultiLin(E_{l_1}^{k_{l_1}},\ldots,E_{l_K}^{k_{l_K}},Y)$$
be such that
$$
D^{k_1,\ldots,k_N}f(x_1\ldots,x_N)(U_1,\ldots,U_K)
=D^{\sum_{i=1}^K k_i}f(x_1,\ldots,x_N)
(\iota_{l_1}^{k_{l_1}}(U_1),\ldots,\iota_{l_K}^{k_{l_K}}(U_K))$$
for $(U_1,\ldots,U_K)\in E_{l_1}^{k_{l_1}},\ldots,E_{l_K}^{k_{l_K}}$.
We also let
$$\chi^{k_1,\ldots,k_N}(f)=\sup\{||D^{k_1,\ldots,k_N}f(x_1,\ldots,x_N)||:
(x_1,\ldots,x_N)\in M_1\times\cdots\times M_N\}.$$

\begin{ex} Supposed $N=3$, $K=2$ and $(k_1,k_2,k_3)=(3,0,2)$ and
  $$f:E_1\times E_2\times E_3\rightarrow Y.$$
  Then $(l_1,l_2)=(1,3)$
and, if $U_1\in E_1^3$ and $U_2\in E_3^2$ then
$$D^{3,0,2}f(x_1,x_2,x_3)(U_1,U_2)
=D^5 f(x_2,x_2,x_3)((U_1)_1,(U_1)_2,(U_1)_3,(U_2)_1,(U_2)_2).$$
for $(x_1,x_2,x_3)\in E_1\times E_2\times E_3$.
\end{ex}
e
\subsection{Diffeomorphisms.}
\begin{defn}\label{diffeodefn} Suppose
$$f:X\rightarrow Y.$$
We say {\bf $f$ carries
$X$ diffeomorphically onto $Y$} if
\begin{itemize}
\item[(i)] there are $m\in\N^+$ and Euclidean spaces $E$ and $F$
such that $X\in\cM_m(E)$ and $Y\in\cM_m(F)$;
\item[(ii)] $f$ is smooth and univalent;
\item[(iii)] $\op{rng}{f}=Y$ and $f^{-1}$ is smooth.
\end{itemize}
\end{defn}

\begin{prop}\label{diffeo} Suppose
\begin{itemize}
\item[(i)] $m\in\N^+$ and $E$ and $F$ are Euclidean spaces;
\item[(ii)] $X\in\cM_m(E)$, $Y\in\cM_m(F)$ and
$$f:X\rightarrow Y;$$
\item[(iii)]  $f$ is smooth and univalent;
\item[(iv)] $\op{rank}{Df(a)}=m$ for each $a\in M$;
\item[(v)] $X$ is compact and $X$ and $Y$ are homeomorphic.
\end{itemize}

Then $\op{rng}{f}=Y$ and $f$ carries
$X$ diffeomorphically onto $Y$.
\end{prop}

\begin{rem} This Proposition need not hold if $X$ is not compact.
\end{rem}

\begin{proof} Let $\cC$ be the family of connected components of $Y$.
$\cC$ is finite since $Y$ is compact.

Suppose $a\in X$. Since $f$ is smooth and $\op{rank}{Df(a)}=m$
is elementary that there are $G$ and $H$ such that $G$ is a relatively
open subset of $X$, $a\in G$, $H$ is an open subset of $F$ and
$f[G]=H\cap\op{rng}{f}$. Let $d=\op{dist}{(a,Y\sim H)}$. Then $d>0$
since $Y\sim H$ is compact and $a\not\in Y\sim H$. So there is an open
subset $W$ of $H$ such that $a\in W$ and $W\cap Y=W\cap f[G]$. Thus $\op{rng}{f}$
is open relative to $Y$. So there is a subfamily $\cB$ of $\cC$
such that $\op{rng}{f}=\bigcup\cB.$ Since $\op{rng}{f}$ is homeomorphic to
$X$ and $X$ is homeomophic to $Y$ we find that $\cB=\cC$ so
$\op{rng}{f}=Y$. That $f^{-1}$ is smooth follows from the Inverse Function
Theorem. So $f$ carries $X$ diffeomorphically onto $Y$
\end{proof}

\subsection{$\cL_m$, $\cP_m$ and $\cW_m$.} 
Let $\cL_m=\{L\in\End(\R^n):\op{rank}{L}=m\}$; as is well known,
$\cL_m\in\cM_{m(n-m)}(\End(\R^n))$.

Let
$$\cP_m:\cL_m\rightarrow\mO_m(\R^n)$$
be such that $\cP_m(L)$ is orthogonal projection of $\R^n$ onto $\op{rng}{L}$;
as is well known, $\cP_m$ is smooth and for each $R\in\R^+$ there
is constant $C_{\cP_m,R}$ such that
\begin{equation}\label{CPmR}
||\cP_m(L)-\cP_m(M)||\leq C_{\cP_m,R}||L-M||.
\end{equation}
whenever $L,M\in\cL_m$ and $\max\{||L||,||M||\}\leq R$.

Let
$$\cW_m:\cL_m\rightarrow\R^+$$
be such that $\cW_m(L)=|\altdown{m}L|$ for $L\in\cL_m$.
It is elementary that $\cW_m$ is smooth and that
$$|\cW_m(L)-1|\leq(2^m-1)||L-\cP_m(L)||$$

\subsection{The second fundamental form of $M$.}

We define
$$B:M\rightarrow{\bigodot}^2(\R^n,\R^n)$$
by requiring that, for any $a\in M$,
$$B(a)(D\phi(\psi(a)(s),D\phi(\psi(a)(t))=\nu(a)(D^2\phi(\psi(a))(s,t))$$
where $U,V,\psi,\psi$ are as in Definition \ref{mani} and $s,t\in\R^m$
and that
$$B(a)(u,v)\bullet w=0\quad\hbox{if either $u\in\Nor{M,a}$ or $v\in\Nor{M,a}$
 or $w\in\Tan{M,a}$.}$$
One calls $B$ the {\bf second fundamental form of $M$}.

\begin{prop} Suppose $a\in M$, $u,v\in\Tan{M,a}$. Then
$$D\tau(a)(u)(v)=B(a)(u,v)=-D\nu(a)(u)(v).$$
\end{prop}
\begin{proof} Let $U,V,\phi,\psi$ are as in Definition \ref{mani},
let $b=\psi(a)$, note that $\phi(b)=a$ and let $s,t\in\R^m$
be such that $u=D\phi(b)(s)$ and $v=D\phi(b)(t)$. Then
\begin{equation}\label{eq}
\nu(\phi(x))(D\phi(x)(t))=0\quad\hbox{for $x\in U$.}
\end{equation}
We differentiate (\ref{eq}) at $x=b$ in the direction $s$ and find that
\begin{align*}
0&=D\nu(a)(D\phi(b)(s))(D\phi(b)(t))+\nu(a)(D(D\phi)(b)(s,t))\\
&=D\nu(a)(D\phi(b)(s))(D\phi(b)(t))+\nu(a)(D^2\phi(b)(s,t))\\
&=D\nu(a)(u)(v)+B(a)(u,v).
\end{align*}
\end{proof}

We define
$$A:M\rightarrow\Lin(\R^n,\Lin(\R^n))$$
by requiring that, for any $a\in M$,
$$A(a)(u)(v)\bullet w=u\bullet B(a)(v,w)\quad\hbox{whenever $u,v,w\in\R^n$.}$$
It follows that
$$A(a)(u)(v)\bullet w=A(a)(u)(w)\bullet v\quad\hbox{if $u,v,w\in\R^n$.}$$

We have
$$||B||=\sup\{|B(a)(u,v)\bullet w|:a\in M \ \hbox{and} \ u,v,w\in\Sm\};$$
It is elementary that $||B||<\infty$ if $M$ is compact and that
$||B||=0$ if and only if for each point $a\in M$ there are an open neighborhood
$V$ of $a$ and $S\in\mG_m(\R^n)$ such that $M\cap V=(a+S)\cap V$.
We let
$$R={1\over||B||}\in[0,\infty];$$
it is elementary that $0<R<\infty$ if $M$ is compact.

\section{Vector bundles.}

\begin{thm} $\mO_m(\R^n)\in\cM_{m(n-m)}(\End(\R^n)).$
\end{thm}
\begin{proof} Let
$$\cP(l)=l\circ l-l\quad\hbox{for $l\in\mathbbm{Sym}(\R^n)$.}$$
Suppose $p\in\mO_m(\R^n)$. Suppose $l\in\op{ker}{D\cP(p)}$.
We have
$$0=D\cP(p)(l)=l\circ p+l\circ l-l=l\circ p-p^\perp\circ l;$$
this implies that $p\circ l\circ p=0$ and $p^\perp\circ l\circ p^\perp=0$
so $l=p\circ l\circ p^\perp+p^\perp\circ l\circ p$
and since $l^*=l$ we have $p^\perp\circ l\circ p=(p\circ l\circ p^\perp)^*$.
\end{proof}

Suppose $l,N$ are integers and $1\leq l\leq N$.
Suppose
$$\gamma:M\rightarrow\mO_l(\R^N)$$
and $\gamma$ is smooth. We let
$$\bE(\gamma)=\{(x,y)\in M\times\R^N:y\in\op{rng}{\gamma(x)}\}$$
and for $r\in\R^+$ we let
$$\bE(\gamma,r)=\{(x,y)\in \bE(\gamma):|y|<r\}.$$

We let
$$\Tan{M}=\bE(\tau)
\quad\hbox{and we let}\quad
\Nor{M}=\bE(\nu).$$

\begin{prop} Suppose $a\in M$ and $u\in\Tan{M,a}$. Then
$$D\gamma(a)(u)=(D\gamma(a)(u))^*$$
and
$$D\gamma(a)(u)=\gamma(a)\circ D\gamma(a)(u)\circ\gamma(a)^\perp
+\gamma(a)^\perp\circ D\gamma(a)(u)\circ\gamma(a).$$
\end{prop}
\begin{proof} The first equation follows from $\gamma(x)^*=\gamma(x)$
for all $x\in M$.

Since $\gamma(x)\circ\gamma(x)=\gamma(x)$ for $x\in M$ we have
$$D\gamma(a)(u)\circ\gamma(a)+\gamma(a)\circ D\gamma(a)(u)=D\gamma(a)(u)$$
for $u\in\R^n$.
Premultiplying and postmultiplying this equation by $\gamma(a)$ we
infer that $\gamma(a)\circ D\gamma(a)(u)\circ\gamma(a)=0$.
Premultiplying and postmultiplying this equation by $\gamma(a)^\perp$ we
infer that $\gamma(a)^\perp\circ D\gamma(a)(u)\circ\gamma(a)^\perp=0$.

\end{proof}

\begin{thm}\label{Egamma} $\bE(\gamma)\in\cM_{m+l}(\R^n\times\R^N)).$
Moreover, for  any $(x,y)\in\bE(\gamma)$ we have
\begin{equation}\label{tangamma}
\begin {split}
\Tan{\bE(\gamma),(x,y)}&=\{(u,v)\in\Tan{M,x}\times\R^N:
\gamma(x)^\perp(v)=D\gamma(x)(u)(y)\}\\
&=\{(u,D\gamma(x)(u)(y)):u\in\Tan{M,x}\}\oplus\{(0,v):v\in\op{rng}{\gamma(x)}\}
\end{split}
\end{equation}
\end{thm}
\begin{proof} We have
$$\Tan{M\times\R^N,(x,y)}=\Tan{M,x}\oplus\R^N
\quad\hbox{for $(x,y)\in M\times\R^N$.}$$
Let $F:M\times\R^N\rightarrow\R^N$ be such that
$$F(x,y)=\gamma(x)^\perp(y)\quad\hbox{for $(x,y)\in M\times\R^N$.}$$
Then $\bE(\gamma)=\{(x,y)\in M\times\R^N:F(x,y)=0\}$.

Suppose $(x,y)\in\bE(\gamma)$. We have
$$DF(x,y)(u,0)=-D\gamma(x)(u)(y)\in\op{rng}{\gamma(x)^\perp}\quad\hbox{for $u\in\Tan{M,x}$}$$
as well as
$$DF(x,y)(0,v)=\gamma(x)^\perp(v)\quad\hbox{for $v\in\R^N$};$$
it follows that
$$\op{rng}{DF(x,y)}=\op{rng}{\gamma(x)}^\perp;\quad\quad
\op{dim}{\op{rng}{DF(x,y)}}=N-l;$$
and
$$\op{ker}{DF(x,y)}=\{(u,v)\in\Tan{M,x}\times\R^N:
\gamma(x)^\perp(v)=D\gamma(x)(u)(y)\}.$$
The Theorem now follows from the Implicit Function Theorem.
\end{proof}

\begin{cor} Suppose $Z\in\cE(\gamma)$, $x\in M$ and $u\in\Tan{M,x}$. Then
$$D\gamma(x)(u)(Z(x))=\gamma^\perp(x)(DZ(x)(u)).$$
\end{cor}
\begin{proof} $(u,DZ(x)(u))\in\Tan{\bE(\gamma),(x,Z(x))}$.
\end{proof}

If $Z\in\cE(\gamma)$ we let
$$D_\gamma Z:M\rightarrow \Lin(\R^n,\R^l)$$
be such that $D_\gamma Z(x)=\gamma(x)\circ DZ(x)$ for $x\in\R^n$.
Evidently, $\op{rng}{D_\gamma Z(x)}\subset\op{rng}{\gamma(x)}$ and
$\Nor{M,x}\subset\op{ker}{D_\gamma Z(x)}$. The preceding Corollary implies that
\begin{equation}\label{Dgamma}
DZ(x)=D_\gamma Z(x)+D\gamma(x)(u)(Z(x))
\quad\hbox{for $x\in M$.}
\end{equation}

\begin{cor}\label{TanNor} If $(x,y)\in\Tan{M}$ then
$$
\Tan{\Tan{M},(x,y)}
=\{(u,B(x)(u,y)+v):(u,v)\in\Tan{M,x}\times\Tan{M,x}\}$$
and if $(x,y)\in\Nor{M}$ then
$$\Tan{\Nor{M},(x,y)}
=\{(u,-A(x)(y)(u)+v):u\in\Tan{M,x}\times\Nor{M,x}\}.$$
\end{cor}

\begin{thm}[Weingarten Equations] Suppose $a\in M$ and $u\in\Tan{M,a}$. Then
\begin{equation}\label{WN}
\tau(a)(DN(a)(u))=-A(a)(N(a))(u)\quad\hbox{if $N\in\cE(\nu)$}
\end{equation}
and
\begin{equation}\label{WT}
\nu(a)(DT(a)(u))=B(a)(T(a),u)\quad\hbox{if $T\in\cE(\tau)$.}
\end{equation}
\end{thm}
\begin{proof} Suppose $N\in\cE(\nu)$. Then $\tau(x)(N(x))=0$ for $x\in M$;
differentiating at $x=a$ we find that
$$\tau(a)(DN(a)(u))=-D\tau(a)(u)(N(a))=-A(a)(N(a))(u).$$

Suppose $T\in\cE(\tau)$. Then $\nu(x)(T(x))=0$ for $x\in M$. Differentiating
at $x=a$ we find that
$$\nu(a)(DT(a)(u)=-D\nu(a)(T(a))(u)=D\tau(a)(T(a))(u)=B(a)(T(a),u).$$
\end{proof}

\subsection{The mean curvature normal of $M$.}\label{mean}
We let
$$H\in\cE(\nu)$$
be such that
$$H(x)\bullet u=\op{trace}{A(x)(u)}\quad\hbox{for $x\in M$ and $u\in\Nor{M,x}$.}$$
$H$ is called the {\bf mean curvature normal of $M$}.

\begin{prop}\label{firstvar} Suppose
\begin{itemize}
\item[(i)] $I$ is an open interval, $0\in I$,
$$h:I\times M\rightarrow \R^n$$
and $h$ is smooth;
\item[(ii)] for each $t\in I$, $M_t=h_t[M]\in\cM_m(\R^n)$ and
$h_t$ carries $M$ diffeomorphically onto $M_t$;
\item[(iii)] for each $t\in I$, $X_t:M_t\rightarrow\R^n$ is such that
$$X_t(y)={d\over d\epsilon}h_{t+\epsilon}(h_t^{-1}(y))|_{\epsilon=0}$$
for $y\in M_t$.
\end{itemize}
Then $X_t$ is smooth and
$${d\over dt}\cH^m(h_t[M])=\int X_t\bullet H_t\,d||M_t||$$
for each $t\in I$.
\end{prop}
For the proof see the literature.

\section{Reach.}\label{reach}

Suppose $M\in\cM_m(\R^n)$.

\begin{defn} We let
$$\Xi:\Nor{M}\rightarrow\R^n$$
be such that
$\Xi(x,y)=x+y\quad\hbox{for $(x,y)\in\Nor{M}$}$
and we let
$$\br=\sup\left(\{0\}\cup\{r\in(0,\infty):\Xi|\bE(\nu,r) \ \hbox{is univalent.}\}\right)$$
\end{defn}

\begin{prop}\label{Xismooth} $\Xi$ is smooth.
\end{prop}
\begin{proof} $\Xi$ is the restriction to $\bE(\nu)$ of the smooth map
$$\R^n\times\R^n\ni(x,y)\mapsto x+y\in\R^n.$$
\end{proof}

\begin{prop} Suppose $0<r<\infty$.
The following statements are equivalent:
\begin{itemize}
\item[(i)] $r\leq\br$;
\item[(ii)] $M\cap\ob{a+v}{|v|}=\{a\}$ whenever $(a,v)\in\Nor{M}$
and $0<|v|<r$;
\item[(iii)] $|\nu(a)(x-a)|\leq{|x-a|^2\over 2r}$ whenever $a,x\in M$.
\end{itemize}
\end{prop}
\begin{proof} Suppose (i) holds, $a\in M$ and $u\in\Nor{M,a}\cap\Sm$.
Let $I=\{s\in(0,r):a+su\not\in M$\}. Suppose $I$ were not empty.
Then $0\leq\inf I<r$. We have $\inf I>0$ since $u\not\in\Tan{M,a}$.
Let $t=\inf I/2$. Then $\Xi(a,tu)=a+tu=\Xi(a+2tu,-tu)$ which is impossible.
Suppose $s\in I$ and $x\in M$ is such that 
$\op{dist}{(a+su,M)}=|(a+su)-x|$. Then $(a+su)-x\in\Nor{M,x}$ so
$\Xi(x,a+su)=\Xi(a,su)$ so, as $s\leq r\leq\br$ we have
$M\cap\ob{a+su,s}=\{a\}$ which implies that $\ob{a+ru,r}\cap M=\emptyset$.
Thus (ii) holds.

It is a simple matter to verify that (ii) implies (i)
and that (ii) and (iii) are equivalent.
\end{proof}

\begin{prop} $\br\leq {1\over||B||}$.
\end{prop}
\begin{proof} Suppose $a\in M$. Suppose
$U,V,\phi,\psi$ are as in Definition \ref{mani} and $u\in\R^m$.
We have
$$B(a)(D\phi(0)(u),D\phi(0)(u))=D^2\phi(0)(u,u)$$
as well as, for $t$ sufficiently small,
$$\phi(tu)=\phi(0)+tD\phi(0)(u)+{t^2\over 2}D^2\phi(0)(u,u)+\rho(t)$$
where $\rho(t)/t^2\rightarrow 0$ as $t\rightarrow 0$.
So
\begin{align*}
{\nu(a)(\phi(tu)-\phi(0))\over t^2}
&={1\over 2}\nu(a)(D^2\phi(0)(u,u))+{\nu(a)(\rho(t))\over t^2}\\
&={1\over 2}(B(a)(D\phi(0)(u),D\phi(0)(u)+{\nu(a)(\rho(t))\over t^2}.
\end{align*}
But
$${|\nu(a)(\phi(tu)-\phi(0))|\over t^2}\leq{|\phi(tu)-\phi(0)|^2\over 2r\,t^2}
\rightarrow{|D\phi(0)(u)|^2\over 2r}\quad\hbox{as} \ t\rightarrow 0.$$
Thus $|B(a)(D\phi(0)(u),D\phi(0)(u)|\leq{|D\phi(0)(u)|^2\over r}$.

\end{proof}

\begin{thm} Suppose $(a,b)\in\Nor{M}$ and $(u,v)\in\Tan{M,a}\times\Nor{M,a}$.
Then
\begin{equation}\label{Xi1}
D\Xi(a,b)(u-A(a)(b)(u),v)=u-A(a)(b)(u)+v
\end{equation}
and
\begin{equation}\label{Xi2}
|D\Xi(a,b)(u-A(a)(b)(u),v)|^2=|u-A(a)(b)(u)|^2+|v|^2.
\end{equation}
In particular, if $|b|<1/||B||$, then
\begin{equation}\label{Xi3}
|D\Xi(a,b)(u-A(a)(b)(u),v)|\geq|u-A(a)(b)(u)|\geq 1-||B||\,|v|>0
\end{equation}
and $D\Xi(a,b)$ carries $\Tan{\Nor{M},(a,b)}$ isomorphically onto
$\R^n$.
\end{thm}
\begin{proof} $(u-A(a)(b)(u),v)\in\Tan{\Nor{M,a},(a,b)}$ by Corollary \ref{TanNor}.

Let $F(x,y)=x+y$ for $(x,y)\in\R^n\times\R^n$. Since $\Xi=F|\Nor{M}$,
(\ref{Xi1}) follows from Corollary \ref{TanNor}. (\ref{Xi2}) holds
since $u-A(a)(b)(u)(v)\in\Tan{M,a}$ and $v\in\Nor{M,a}$. If $|b|<1/||B||$
then $|u-A(a)(b)(u)|>1-||B||\,|v|$ and (\ref{Xi3}) follows from (\ref{Xi2}).
(\ref{Xi3}) implies $D\Xi(a,b)|\Tan{\Nor{M},(a,b)}$ is univalent.
This and the fact that $\Nor{M}\in\cM_n(\R^n\times\R^n)$ implies the
final conclusion.
\end{proof}

\begin{defn} For $r\in(0,\infty]$ we set
$$N(r)=\{z\in\R^n:\op{dist}{(z,M)}<r\}.$$
\end{defn}
Obviously, $N(r)$ is an open subset of $\R^n$.

\begin{thm}\label{Xidiffeo} Suppose $\br>0$. Then $\Xi$ carries $\bE(\nu,\br)$
diffeomorphically onto $N(\br)$.
\end{thm}
\begin{proof} $\Xi|\bE(\nu,\br)$ is univalent by the definition of $\br$.
  If $(a,b)\in\bE(\nu,\br)$ then $D\Xi(a,b)$ carries $\Tan{\Nor{M},(a,b)}$
  isomorphically onto $\Tan{N(\br),a+b}=\R^n$ by the preceding Theorem.
\end{proof}

\begin{prop} $\br=\infty$ if and only if $M$ is a subset
of some member of $\mG_m(\R^n)$.
\end{prop}
\begin{proof} Suppose $\br=\infty$ and $a\in M$. It follows from (ii) of the
preceding Proposition that $M\subset\Tan{M,a}$.
\end{proof}

\begin{prop} Suppose $M\subset\Sm$. If $a\in M$ then
$A(a)(u)=u$ for $u\in\Tan{M,a}$; in particular, $||B||\geq 1$
and $\br\leq 1$.
\end{prop}
\begin{proof} Let $X(x)=x$ for $x\in M$. Then $X\in\cE(\nu)$
and, if $u\in\Tan{M,a}$, $DX(a)(u)=u$ so
so $A(a)(u)=u$ by a (\ref{WN}); also, if $|u|=1$, then
$||B||\geq||A(a)(u)||= 1$.
\end{proof}

\begin{lem} \label{reachlemma} Suppose
\begin{itemize}
\item[(i)] $Y_i$, $i=1,2$, is a Euclidean space;
\item[(ii)] $N_i\in\cM_m(Y_i)$ for $i=1,2$;
\item[(iii)] $K$ is a compact subset of $N_1$;
$F:N_1\rightarrow N_2$ is smooth; $F|K$ is univalent and,
for each $x\in K$,
$DF(x)$ carries $\Tan{N_1,x}$ isomorphically onto $\Tan{N_2,F(x)}$;
\item[(iv)] $U_\epsilon$, $0<\epsilon<\infty$, is an open subset of $N_1$
such that
\begin{itemize}
\item[(a)] if $0<\epsilon_1<\epsilon_2<\infty$ then $U_{\epsilon_1}\subset
U_{\epsilon_2}$;
\item[(b)] if $V$ is an open subset of $N_1$ and $K\subset V$ then
$U_\epsilon\subset V$ for some $\epsilon>0$;
\item[(c)] $K=\bigcap_{\epsilon>0}U_\epsilon$.
\end{itemize}
\end{itemize}

Then
\begin{itemize}
\item[(v)] there is $\epsilon>0$ such that $DF(x)$ carries $\Tan{N_1,x}$ isomorphically onto $\Tan{N_2,F(x)}$
for $x\in U_\epsilon$ and such that $F|U_\epsilon$ is univalent;
\item[(vi)] if $\epsilon$ is as in (v) then $F[U_\epsilon]$ is an open neighborhood in $N_2$ of $F[K]$ and $F$ carries $U_\epsilon$ diffeomorphically onto
$F[U_\epsilon]$.
\end{itemize}
\end{lem}
\begin{rem} What follows is a paraphrase of the argument on page 116 of
\cite{MS}.
\end{rem}
\begin{proof} Let $V_1$ be the set of $x\in N_1$ such that $DF(x)$ carries
$\Tan{N_1,x}$ isomorphically onto $\Tan{N_2,F(x)}$.
Then $V_1$ is an open subset
of $N_1$ so there is by (iii)(b) $\epsilon>0$ such that $U_\epsilon\subset V_1$.

Suppose (v) is false. Let $N_*$ be a positive integer such that the
closure of $K$ of $U_{1/N_*}$ in $N_1$ is compact.
For each integer $N$ such that $N\geq N_*$ there are $a_N,b_N\in U_{1/N}$ such that
$a_N\neq b_N$ and $F(a_N)=F(b_N)$. Since $K$ is compact there are
 $a_*$ and $b_*$ in $K$ which are accumulation points
of $\{a_N:N=1,2,\ldots\}$ and $\{b_N:N=1,2,\ldots\}$, respectively.
Then $F(a_*)=F(b_*)$ which implies $\lim_{N\rightarrow\infty}|a_N-b_N|=0$.
Thus $F$ would not be univalent in any open neighborhood
of $a_*$. But as $DF(a_*)$ carries $\Tan{N_1,a_*}$ isomorphically onto
$\Tan{N_2,F(a_*)}$, $F$ must be univalent in some open neighborhood of $a_*$.

(vi) is an elementary consequence of (v) and the smoothness of $F$. 
\end{proof}

\begin{cor} Suppose $M\in\cM_m(\R^n)$ and $M$ is compact. Then
$\br>0$.
\end{cor}
\begin{proof} Let
$$Y_1=\R^n\times\R^n;\quad Y_2=\R^n;\quad N_1=\Nor{M};\quad N_2=\R^n;$$
let $F=\Xi$ and,
for $\epsilon>0$, let $U_\epsilon=\bE(\gamma,\epsilon)$; and let
$K=\{(x,0):x\in M\}\subset N_1$. Then (i)-(iv) of Lemma \ref{reachlemma} hold.
So if $\epsilon$ is as in Lemma \ref{reachlemma} then $\reach{M}\geq\epsilon$.
\end{proof}

\subsection{}
Suppose $M$ is compact.
Since $\Xi$ carries $\bE(\nu,\br)$ diffeomorphically onto $N(\br)$
we may define smooth functions
$$\xi:N(\br)\rightarrow \R^n\quad\hbox{and}\quad\sigma:N(\br)\rightarrow\R^n$$
by requiring that
$$(\xi,\sigma)=(\Xi|\bE(\nu,\br)^{-1}:N(\br)\rightarrow\bE(\nu,\br).$$
Note that
\begin{equation}\label{xisigma}
z=\Xi(\xi(z),\sigma(z))=\xi(z)+\sigma(z)
\quad\hbox{for $z\in N(\br)$.}
\end{equation}

\begin{thm}\label{Dxisigma} Suppose $b\in N(\br)$, $a=\xi(b)$ and $c=\sigma(b)$.
Then
\begin{itemize}
\item[(i)] $a\in M$, $c\in\Nor{M,a}$, $b=a+c$ and $|c|=\op{dist}{(b,M)}$;
\item[(ii)] $D\xi(b)(u-A(a)(c)(u))=u$ for $u\in\Tan{M,a}$ and
$\op{ker}{D\xi(b)}=\Nor{M,a}$;
\item[(iii)] $D\sigma(b)(u-A(a)(c)(u))=-A(a)(c)(u)$ for $u\in\Tan{M,a}$ and
$D\sigma(b)(v)=v$ for $v\in\Nor{M,b}$;
\item[(iv)] $||D\xi(b)-\tau(a)||\leq{||B||\,|c|\over 1-||B||\,|c|}\leq{\br|c|\over \br-|c|}$.
\end{itemize}
\end{thm}
\begin{rem} Keep in mind that
$$\Tan{M,a}\ni u\mapsto u-A(a)(c)(u)\in\Tan{M,a}$$
is a linear isomorphism since $||B||\,|c|<1$.
\end{rem}
\begin{proof} We have $\xi\circ\Xi(x,y)=x$ and $\sigma\circ\Xi(x,y)=y$
for $(x,y)\in\bE(\nu,\br)$. (i)-(iii) now follow from (\ref{Xi1})
and the Chain Rule.

Let $l\in\End(\Tan{M,a})$ and $L\in\Lin(\Tan{M,a},\R^n)$ be such that $l(u)=A(a)(c)(u)$
and $L(u)=u-l(u)$ for $u\in\Tan{M,a}$. Then $||l||\leq ||B||\,|c|<1$
and $||L||\geq 1-||l||$ so $L$ is invertible and $||L^-1||\leq(1-||l||)^{-1}$.
Let $i$ be the identity map of $\Tan{M,a}$. (ii) implies that
$D\xi(b)\circ L=i$ so $D\xi_b|\Tan{M,a}=L^{-1}$ which implies that
$$||D\xi(b)|\Tan{M,a}-i||
=||L^{-1}-i||\leq||L^{-1}||\,||i-L||\leq{||l||\over 1-||l||}
\leq{||B||\,|c|\over 1-||B||\,|c|}.$$
(iv) follows since (ii) and (iii) imply that $\Nor{M,a}\subset\op{ker}{D\xi(b)}$.
\end{proof}

\begin{thm}\label{Dxi} Suppose $z\in N(\br)$. Then
\begin{equation}\label{xitau}
D\xi(\xi(z))=\tau(\xi(z))
\end{equation}
and
\begin{equation}\label{Dxitau}
||D\xi(z)-\tau(\xi(z))||\leq\chi^2(\xi)\op{dist}{(z,M)}.
\end{equation}
\end{thm}
\begin{proof} (\ref{xitau}) follows from (ii) of Theorem \ref{Dxisigma}.
Let $C(\epsilon)=(1-\epsilon)z+\epsilon\xi(z)$ for $\epsilon\in(0,1)$. 
Since
$\ob{z}{\dist{z,M}}\subset N(r)=\op{dmn}{\xi}$ we find that
$\op{rng}{C}\subset\op{dmn}{\xi}$ so that
$$D\xi(z)-\tau(\xi(z))
=D\xi(z)-D\xi(\xi(z))
=\int_0^1 D(D\xi)(C(\epsilon))(C^\prime(\epsilon))\,d\epsilon$$
so, as $C^\prime(\epsilon)=z-\xi(z)$ for $\epsilon\in(0,1)$, (\ref{Dxitau})
holds.
\end{proof}

\begin{thm}\label{Lip} Suppose $f$ is a smooth function on $M$ with values in a
normed space $Y$. Then
$$\op{Lip}{f}\leq \left(1+{1\over\br}\right)\chi^{[1]}(f).$$
\end{thm}
\begin{proof}  Suppose $0<\lambda<1$ and
$|a-b|\leq 2\lambda\br$. Let $C(t)=(1-t)a+tb$ for $t\in(0,1)$.
Then $\op{dist}{(C(t),M)}<\lambda\br$ for $t\in(0,1)$ so $\op{rng}{C}\subset\op{dmn}{\xi}$ and, by Theorem \ref{Dxisigma}(ii),
$$||D\xi(C(t))||
\leq {1\over 1-||B||\op{dist}{(C(t),M)}}
\leq{1\over 1-\lambda}$$
so
$${d\over dt}(f\circ\xi\circ C)(t)
\leq \chi^1(f)||D\xi(C(t))||\,|C^\prime(t)|
\leq {\chi^1(f)\over 1-\lambda}|a-b|.$$
Integrating this inequality from $0$ to $1$ we infer that
$$|f(a)-f(b)|\leq {\chi^1(f)\over 1-\lambda}|a-b|.$$

On the other hand, if $|a-b|>2\lambda\br$ then
$$|f(a)-f(b)|\leq \op{diam}{\op{rng}{f}}
\leq 2\chi^0(f)
<{\chi^0(f)\over \lambda\br}|a-b|.$$

Thus
$$|f(a)-f(b)|\leq\chi^{[1]}\max\left\{{1\over 1-\lambda},{1\over\lambda\br}\right\}|a-b|$$
so we should set $\lambda=1/(1+\br)$.
\end{proof}

\subsection{$\theta$ and $\iota$.} Let
$$\theta:\Sm\rightarrow\mO_{n-1}(\R^n)
\quad\hbox{and}\quad
\iota:\Sm\rightarrow\mO_1(\R^n)$$
be the tangent map of $\Sm$ and normal map of $\Sm$, respectively.

\begin{cor} Suppose $(x,v)\in\bE(\theta)$. Then
$$D\iota(x)(u)=$$
$$\Tan{\bE(\theta),(x,v)}
=\{(u,w)\in\Tan{\Sm,x}\times\R^n:w\bullet x=-u\bullet v\}.$$
\end{cor}

\section{$\bbeta$ and $\bo$.}
We define
$$\bbeta:\R^n\sim\{0\}\rightarrow\Sm
\quad\hbox{and}\quad
\bo:\R^n\sim\{0\}\rightarrow\mathbb{G}_{n-1}(\R^n)$$
at $x\in\R^n\sim\{0\}$
by letting $\bbeta(x)=|x|^{-1}x$
and letting $\bo(x)$ be orthogonal projection of $\R^n$
onto $\Tan{\Sm,x}=\{w\in\R^n:w\bullet x=0\}$.

\begin{prop}\label{Lipbbeta} Suppose $x,y\in\Sm\sim\{0\}$. Then
\begin{equation}\label{eta}
|\bbeta(x)-\bbeta(y)|\leq {|x-y|\over\min\{|x|,|y|\}}.
\end{equation}
\end{prop}
\begin{proof} Let $s,t\in\R^+$ and $u,v\in\Sm$ be such that
$x=su$ and $y=tv$. Since both sides of (\ref{eta}) are symmetric in $x$ and $y$
we may suppose that $0<s=|x|\leq|y|=t$. Then (\ref{eta}) amounts to
$$|u-v|\leq{|su-tv|\over s}=|u-\lambda v|$$
where we have set $\lambda=t/s\geq 1$. We have
$$|u-\lambda v|-|u-v|
=1-2\lambda u\bullet v+\lambda^2-(1-2 u\bullet v+1)
=(\lambda-1)(\lambda+1-2 u\bullet v)
\geq 0.$$
\end{proof}

\begin{prop}\label{Dbbeta} For $x\in\R^n\sim\{0\}$ we have
$$D\bbeta(x)={\bo(x)\over |x|};$$
in particular,
  $$D\bbeta(x)(w)=|x|^{-3}(|x|^2 w-(w\bullet x)x)
  \quad\hbox{for $w\in\R^n$.}$$
\end{prop}
\begin{proof} Suppose $w\in\R^n$. Then
$$D\bbeta(x)(w)={d\over d\epsilon}{x+\epsilon w\over |x+\epsilon w|}\big|_{\epsilon=0}
={w\over|x|}-{1\over|x|^2}\left(w\bullet{x\over|x|}\right)x
=|x|^{-3}(|x|^2 w-(w\bullet x)x).$$
\end{proof}

\begin{prop}\label{bbetaxv} Suppose $x\in\Sm$ and $v\in\R^n$ and $v\bullet x=0$. Then
  $$|\bbeta(x+v)-x|\leq|v|.$$
\end{prop}
\begin{proof} We have
$$|x-\bbeta(x+v)|^2=2\left(1-x\bullet\left(x+v\over|x+v|\right)\right)
=2\left(1-{1\over\sqrt{1+|v|^2}}\right)
\leq|v|^2.$$
\end{proof}

\section{ $f$.}\label{f}

We let
$$f:\bE(\nu,\br)\rightarrow\R^n$$
be such that
$$f(x,y)=(1+y\bullet x)\bbeta(x+y)\quad\hbox{for $(x,y)\in\bE(\nu,\br)$.}$$
We let
$$\upsilon:N(\br)\rightarrow\R$$
be such that $\upsilon(z)=z\bullet\xi(z)$ for $z\in N(\br)$.

\begin{rem} On page 221 of \cite{AA81} $f$ there is based on the exponential map
of $\Sm$ applied to $\Xi|\bE(\gamma)$. The $f$ here is as in [FN]
which I believe is easier to work with. The modifications to required to
accomodate this $f$ in \cite{AA81} are straightforward.
\end{rem}

\begin{prop} $\upsilon>0$.
\end{prop}
\begin{proof} Suppose $z\in N(\br)$ and $x=\xi(z)$. Then
$$z\bullet\xi(z)=(x+(z-x))\bullet x=1+(z-x)\bullet x
\geq 1-|z-x|>1-\br\geq 0.$$
\end{proof}

We define the smooth function
$$\rho:N(\br)\rightarrow\R^n$$
by setting
$$\rho(z)={z\over\upsilon(z)}-\xi(z)\quad\hbox{for $z\in N(\br)$.}$$

\begin{thm}\label{fdiffeo} $f$ is a diffeomorphism with range an open subset
$V$ of $N(\br)$ and
with inverse $(\xi,\rho)|V$. In particular,
\begin{equation}\label{DxifP}
D\xi(f(a,b))(Df(a,b)(u-A(a)(b)(u),v)=DP(a,b)(u-A(a)(b)(u),v)=u-A(a)(b)(u)
\end{equation}
and
\begin{equation}\label{DrhofQ}
D\rho(f(a,b))(Df(a,b)(u-A(a)(b)(u),v)=DQ(a,b)(u-A(a)(b)(u),v)=v
\end{equation}
for any $(a,b)\in\bE(\nu,\br)$ and $(u,v)\in\Tan{M,a}\times\Nor{M,a}$
and where
$$P(a,b)=a\quad\hbox{and}\quad Q(a,b)=b\quad\hbox{for $(a,b)\in\bE(\nu,\br)$.}$$
\end{thm}
\begin{proof} Suppose $(x,y)\in\bE(\nu,\br)$ and $z=f(x,y)$. Then
\begin{itemize}
\item[(i)] $1+y\bullet x>1-\br\geq 0$ and $|z|=1+y\bullet x$;
\item[(ii)] $z-x\in\Nor{M,x}$;
\item[(iii)] $z\in\op{dmn}{\xi}$ and $x=\xi(z)=x$;
\item[(iv)] $y=\rho(z)$.
\end{itemize}
$y\bullet x\leq|y|<\br$ and $|z|=1+y\bullet x$ so (i) holds.

$\{x,y\}\subset\Nor{M,x}$ so (ii) holds. 

Let $\alpha(r)=(1+r)^2-2{(1+r)^2\over 1-r}+1$ for $r\in[0,1)$. Then
$(1-r)(r^2-\alpha(r))=2(2|y|^2+|y|+1>0$ so $\alpha(|y|)<|y|^2<\br^2$.
We infer from
$1+y\bullet x\leq 1+|y|$ and $|x-y|>1-|y|$ that
\begin{align*}
|z-x|^2
&=\left|(1+y\bullet x)^2-2(1+y\bullet x){x+y\over|x+y|}\bullet x+1\right|^2\\
&=(1+|y|)^2-2{(1+|y|)^2\over 1-|y|}+1\\
&=\alpha(|y|)\\
&<\br^2
\end{align*}
so that $\xi(z)=x$ and (iii) holds.

(i) implies that
$$z=|z|{x+y\over|x+y|};$$
Taking the inner product with $\xi(z)=x$ we find that
$$z\bullet\xi(z)={1+y\bullet x\over|x+y|}={|z|\over|x+y|}$$
so that $|x+y|={|z|\over z\bullet\xi(z)}$ and
$$z=|z|{z\bullet\xi(z)\over|z|}(x+y)=(z\bullet\xi(z))(\xi(z)+y)$$
so (iv) holds.

We have shown that $(\xi,\rho)\circ f(x,y)=(x,y)$ for $(x,y)\in\bE(\nu,\br)$.
Thus $f$ is univalent and, by the Chain Rule, $Df(x,y)$
carries $\Tan{\bE(\nu,\br),(x,y)}$ isomorphically onto $\R^n$ so
$f$ is an open mapping. The Theorem follows.
\end{proof}

\subsection{$\gamma$.}\label{gamma} Let
$$\gamma(x)=\nu(x)-\iota(x)\quad\hbox{for $x\in M$;}$$
evidently, $\gamma:M\rightarrow\mO_{n-m-1}(\R^n)$, $\gamma$ is smooth
and $\op{rng}{\gamma(x)}=\Nor{M,x}\cap\Tan{\Sm,x}$ for any $x\in M$.

\begin{prop} Suppose $a\in M$, $u\in\Tan{M,a}$ and $y\in\Nor{M,a}$. Then
\begin{equation}\label{auy}
D\gamma(a)(u)(y)=
-A(a)(y)(u)-
\begin{cases}
u&\hbox{if $y=a$,}\cr
0&\hbox{if $y\in\op{rng}{\gamma(a)}$}
\end{cases}
\end{equation}
and
\begin{equation}\label{DgammaW}
DW(a)=D_\gamma W(a)-A(a)(W(a))\quad\hbox{for $W\in\cE(\gamma)$.}
\end{equation}

Moreover,
\begin{equation}\label{auyTan}
\Tan{\bE(\gamma),(a,v)}
=\{(u,A(a)(a)(u):u\in\Tan{M,a}\}
\oplus\{(0,w):w\in\op{rng}{\gamma(a)}\}.
\end{equation}
\end{prop}
\begin{proof} We have
$$\gamma(x)(y)=\nu(x)(y)-(y\bullet x)x\quad\hbox{for $x\in M$
and $y\in\R^n$.}$$
For $u\in\Tan{M,a}$ we have
$$D\gamma(a)(u)(y)=-A(a)(y)(u)-(y\bullet u)a-(y\bullet a)u.$$
Since 
$$(y\bullet u)a+(y\bullet a)u
=\begin{cases}
u&\hbox{if $y=a$,}\cr
0&\hbox{if $y\in\op{rng}{\gamma(a)}$}
\end{cases}$$
(\ref{Dgamma}) implies (\ref{auy}). (\ref{auyTan}) now follows from (\ref{tangamma}).
\end{proof}

The following Theorem is a Corollary of Theorem \ref{fdiffeo}.
\begin{thm} $f$ carries $\bE(\gamma,\br)$ diffeomorphically
onto $V\cap\Sm$.
\end{thm}

\begin{prop}\label{Drho} Suppose $x\in M$. Then
$$D\rho(x)(w)=\gamma(x)(w)\quad\hbox{whenever $w\in\Nor{M,x}\cap\Tan{\Sm,x}$.}$$
\end{prop}
\begin{proof} We have $D\xi(x)=\tau(x)$.

Suppose $w\in\Nor{M,x}\cap\Tan{\Sm,x}$. Then
$$D\upsilon(x)(w)=w\bullet x+x\bullet D\xi(x)(w)=0.$$
so, as $\upsilon(x)=1$,
$$D\rho(x)(w)=w-D\xi(x)(w)=\nu(x)(w)=\gamma(x)(w).$$
\end{proof}

\subsection{$F$.}\label{F}

Whenever $W\in\cE(\gamma,\br)$ we define
$$F(W):M\rightarrow\Sm$$
by letting $F(W)=f\circ(i_M,W)$; it follows from Theorem \ref{fdiffeo} that
$F(W)$ carries $M$ diffeomorphically onto its range; in particular,
\begin{equation}\label{Fprop}
\op{rng}{F(W)}\in\cM_m(\R^n),\quad\xi\circ F=i_M \quad\hbox{and}\quad \rho\circ F(W)=W.
\end{equation}

\section{Families of submanifolds of $\Sm$.}

We suppose throughout this section that $\cJ$ is a Euclidean space.
For each $\epsilon\in\R^+$ and $s\in\cJ$ we let
$$\cJ(\epsilon)=\{t\in\cJ:|t|<\epsilon\}
\quad\hbox{and we let}\quad
\cJ(s,\epsilon)=\{t\in\cJ:|s-t|<\epsilon\}.$$

We suppose
\begin{itemize}
\item[(i)] $\cS$ is a open subset of $\cJ$ and $0\in\cS$;
\item[(II)] $M\in\cM_m(\R^n)$, $M\subset\Sm$ and $M$ is compact;
\item[(III)] $\bz:\cS\times M\rightarrow\Sm$ and $\bz$ is smooth;
\item[(IV)] $\bz_s$ is a diffeomorphism for each $s\in \cS$;
\item[(IV)] $\bz_0(p)=p$ for $p\in M$.
\end{itemize}

Let
$$\br$$
be the reach of $M$.

Let
$$M_s=\bz_s[M]\quad\hbox{for each $s\in\cS$;}$$
it follows that $M_s\in\cM_m(\R^n)$, $M_s$ is compact
and $M_s\subset\Sm$.

We let
$$T:\cS\times M\rightarrow\cL_m$$
be such that $T(s,p)=\cP_m(D\bz_s(p))$ for $(s,p)\in\cS\times M$.

For each $s\in \cS$ we let
$$\tau_s, \ \nu_s, \ B_s, \ A_s \quad\hbox{be as in Section \ref{subman};}$$
$$ \br_s, \ \Xi_s, \ N_s(r), \ r\in\R^+, \ \xi_s, \ \sigma_s\quad\hbox{be as in Section \ref{reach}};$$
$$\upsilon_s, \ \rho_s, \ f_s, \ V_s, F_s
 \quad\hbox{be as in Section \ref{f}}$$
with $M$ replaced by $M_s$.

We have
$$T(s,p)=\tau_s(\bz_s(p))\quad\hbox{for $(s,p)\in\bS\times M$.}$$

\subsection{$\epsilon_1$.}

We fix
$$\epsilon_1\in\R^+$$
such that the closure of $\cJ(\epsilon_1)$ is a compact subset of $\cS$.

For each $(j,k)\in\N\times\N$ we let
\begin{equation}\label{bound}
C_{j,k}=\chi^{j,k}(\bz|(\cJ(\epsilon_1)\times M))
\quad\hbox{for any $(j,k)\in\N$.}
\end{equation}
and note that $C_{j,k}<\infty$.

\begin{thm} There is a constant $C_T$ such that
\begin{equation}\label{CT}
||T(s,p)-T(t,q)||\leq C_T(\epsilon)(|s-t|+|p-q|)
\quad\hbox{whenever $s,t\in\cJ(\epsilon_1)$ and $p,q\in M$}
\end{equation}
\end{thm}
\begin{proof} Suppose  $s,t\in\cJ(\epsilon_1)$ and $p,q\in M$. By (\ref{bound})
we have
$$T(s,p)-T(t,p)
\int_0^1{d\over d\lambda}T((1-\lambda)t+\lambda s,p)\,d\lambda
=\int_0^1 DT((1-\lambda)t+\lambda s,p)(s-t,0)\,d\lambda
$$
and, by Theorem \ref{Lip},
$$||T(t,p)-T(t,q)||
\leq\left(1+{1\over\br}\right)\chi^{[1]}(T_t)|p-q|.$$
\end{proof}

Let
$$\ubr=\inf\{\br_s:s\in \cJ(\epsilon_1)\}
\quad\hbox{and let}\quad
\ubb=\sup\{||B_s||:s\in\cJ(\epsilon_1)\}.$$
\begin{thm} We have
$$0<\ubr\leq 1\quad\hbox{and}\quad \ubb\leq{1\over\ubr}.$$
\end{thm}\label{ubr}
\begin{proof} Let $\cH$ be the closure of $\cJ(\epsilon_1)$.
Let
$$Y_1=(\cJ\oplus\R^n)\oplus\R^n,\quad
Y_2=\cJ\oplus\R^n,\quad
N_1=\bE(T),\quad
N_2=\cJ\oplus\R^n,$$
let $F:N_1\rightarrow N_2$
be such that $F((s,p),v)=\bz(s,p)+v$ for $((s,p),v)\in\bE(T)$.
Let $K=\{((s,p),0):(s,p)\in\cH\times M\}$ and for each $\epsilon\in\R^+$
let
$$U_\epsilon=\{((s,p),v)\in\bE(T):\op{dist}{(s,K)}<\epsilon \ \hbox{and}
 \ |v|<\epsilon\}.$$
 Keeping in mind Theorem \ref{Egamma} we find that (i)-(iv)
 of Lemma \ref{reachlemma} hold. Let $\epsilon_*$ be as in (v) of
 Lemma \ref{reachlemma}. Then
 $\epsilon_*\leq\br_s$ for any $s\in\cJ(\epsilon_1)$.

$\ubb\leq 1/\ubr$ since $||B_s||\br_s\leq 1$ for $s\in\cJ$.
\end{proof}


For each $(s,t)\in\cJ(\epsilon_1)$ we let
$$\bc_{s,t}:[0,1]\rightarrow \cJ(\epsilon_1)$$
be such that
$$\bc_{s,t}(\lambda)=t+\lambda(s-t)\hbox{for $\lambda\in[0,1]$.}$$

\begin{prop}\label{zdiff} Suppose $s,t\in \cJ(\epsilon_1)$. Then
$$|\bz_s(p)-\bz_t(p)|\leq C_{1,0}|s-t|\quad\hbox{for any $p\in M$.}$$
Moreover,
$$M_t\subset N_s(r)\quad\hbox{if $C_{1,0}|s-t|<r\leq\ubr$.}$$
\end{prop}
\begin{proof} For any $p\in M$ we have
$$|\bz_s(p)-\bz_t(p)|
\left|\int_0^1 {d\over d\lambda}\bz(\bc_{s,t}(\lambda),p)\,d\lambda\right|
=\left|\int_0^1 D\bz(\bc_{s,t}(\lambda),p)(s-t,0)\,d\lambda\right|
\leq C_{1,0}|s-t|.$$

Suppose $z\in M_t$. Let $p\in M$ be such that $\bz_t(p)=z$. Then
$$\op{dist}{(z,M_s)}\leq|\bz_t(p)-\bz_s(p)|\leq C_{1,0}|s-t|.$$
\end{proof}

Let
$$\delta_1\in\R^+$$
be such that $C_{1,0}\delta_1<\ubr$.

\begin{cor}\label{delta1} If $s,t\in\cJ(\epsilon_1)$
and $|s-t|<\delta_1$ then
$$M_t\subset N_s(\ubr).$$
\end{cor}


\newcommand{\mmu}{\langle \mu \rangle}
We define
$$\mu:\cS\times M\rightarrow\R^+$$
by letting 
$$\mu(s,p)=\sup\{\lambda\in[0,\infty):
|D\bz_s(p)(u)|\geq\lambda|u| \ \hbox{for} \ u\in\Sm\cap\Tan{M,p}\}.$$
We let
$$\mmu=\inf\{\mu(s,p):(s,p)\in\cJ(\epsilon_1)\times M\}$$

We leave as an elementary exercise for the reader to prove the following
Proposition.
\begin{prop} $\mu$ is positive and continuous and $\mmu$ is positive.
\end{prop}

\begin{prop}\label{Cmu} There is a constant
$$C_\mu$$
such that if $s,t\in\cJ(\epsilon_1)$ and $|s-t|<\delta_1$ then
$M_t\subset N_s(\ubr)$ and
$$||T(t,p)-\tau_s(\xi_s(\bz_t(p)))||\leq C_\mu|s-t|
\quad\hbox{for $p\in M$.}$$
\end{prop}
\begin{proof} Suppose $s,t\in\cJ(\epsilon_1)\cap\cJ(\delta_1)$.
From Corollary \ref{delta1}
we infer that
$$M_t\subset N_s(\ubr)\subset\op{dmn}{\xi_s}.$$

Suppose $p\in M$. Then
$$||D\bz_s^{-1}(p)||\leq{1\over\mu(s,p)}\leq{1\over\mmu}.$$
Let $q\in M$ be such that $\xi_s(\bz_t(p))=\bz_s(q)$.
If
$$c=\left(1+{1\over \ubr}\right)\left({1\over\ul{\mu}}\right)$$
we infer from Theorem \ref{Lip} that
\begin{align*}
|p-q|
&=|\bz_s^{-1}(\bz_s(p))-\bz_s^{-1}(\bz_s(q))|\\
&\leq\Lip{\bz_s^{-1}}|\bz_s(p)-\bz_s(q)|\\
&\leq\left(1+{1\over \br_s}\right)\left({1\over\ul{\mu}}\right)|\bz_s(p)-\bz_s(q)|\\
&\leq c|\bz_s(p)-\bz_s(q)|\\
&\leq c(|\bz_s(p)-\bz_t(p)|+|\bz_t(p)-\xi_s(\bz_t(p))|\\
&\leq 2c|\bz_s(p)-\bz_t(p)|\\
&\leq 2c C_{1,0}|s-t|.
\end{align*}
Thus
\begin{align*}
||T(t,p)-\tau_s(\xi_s(\bz_t(p)))||
&=||T(t,p)-T(s,q)||\\
&\leq||T(s,p)-T(t,p)||+||T(t,p)-T(t,q)||\\
&\leq C_{1,0}\chi^{1}(T,\cJ(\epsilon_1)\times M)|s-t|+C_T|p-q|\\
&\leq C_\mu|s-t|
\end{align*}
where
$$C_\mu=\chi^{1}(T,\cJ(\epsilon_1)\times M)|s-t|
+C_T(2c C_{1,0}).$$

\end{proof}

\begin{thm}\label{Lambda} Suppose $s\in\cJ(\epsilon_1)$.
There is $\delta\in(0,\epsilon_1-|s|]$
such that if $t\in\cJ(s,\delta)$ then $t\in\cJ(\epsilon_1)$,
$M_t\subset N_s(\ubr)$ and
$$\Lambda_t \ \hbox{carries $M_t$ diffeomorphically onto $M_s$}$$
where
$$\Lambda:\cJ(s,\delta)\times M\rightarrow M$$
is such that
$$\Lambda(t,q)=\bz_s^{-1}\circ\xi_s\circ\bz_t(q)
\quad\hbox{ for $(t,q)\in\cJ(s,\delta)\times M$.}$$
\end{thm}
\begin{proof} That $s,t\in\cJ(\epsilon_1)$ and $|s-t|<\delta_1$.
Then $M_t\subset N_s(\ubr)$ by Corollary \ref{delta1}.

For each $(t,x)\in\cJ(s,\delta_1)\times M_t$ we let
$\cN(t,x)=\{z\in M_t:\xi_s(z)=x\}$ and we let
$N(t,x)$ equal the cardinality of $\cN(t,x)$.

We choose $\delta_2\in(0,\delta_1]$ such that $C_\mu\delta_2<1$
where $C_\mu$ is as in Proposition \ref{Cmu}. Suppose $p\in M$. By Proposition
\ref{Cmu} we have
$$||T(t,p)-\tau_s(\xi_s(\bz_t(p))||<1
\quad\hbox{ if $t\in\cJ(s,\delta_2)$.}$$
It follows from (ii) of Theorem \ref{Dxisigma}
that $\Tan{M_t,\bz_t(p)}$
meets $\op{ker}{D\xi_s(\bz_t(p))}$ transversely so that
$$\hbox{$D\xi_s(\bz_t(p))$
carries $\Tan{M_t,\bz_t(p)}$ isomorphically onto $\Tan{M_s,\xi_s(\bz_t)(q)}$.}$$
This implies that there is an open neighborhood of $\bz_t(p)$ such that
$\bz_s$ carries $U\cap M_t$ diffeomorphically onto an open neighborhood of
$\bz_s(\bz_t(p))$ in $M_s$. 
Since $M_t$ is compact we find that if $x\in\op{rng}{\xi_s|M_t}$ then
$N_t(x)$ is finite and there are an open neighborhood $U$ of $x$
in $M_s$ and  a map
$\cU$ with domain $\cN_t(x)$ such that $\cU(z)$ is an open neighborhood of
$z$ in $M_t$ which $\xi_s$ carries diffeomorphically onto $U$.
It follows that $\cN(t,x)$ is finite and that $N_t$ is locally constant.

We choose $\delta_3\in(0,\delta_2]$ as follows. Let $\cC$ be the family of
connected components of $M_s$; since $M_s$ is compact $\cC$ is finite.
For $x\in M_s$ we let $C_x\in\cC$ be such that $x\in C_x$.
If $\cC$ has one member we let $\delta_2=\delta_1$. Otherwise we let
$$d=\inf\{|x-y|:x,y\in M_s \ \hbox{and} \ C_x\neq C_y\}.$$
We choose $\delta_3\in(0,\delta_2]$ such that
$C_{1,0}\delta_2<d$.
Suppose $|s-t|<\delta_3$ and $x\in M_s$. There is $z\in M_t$
such that $|x-z|=\op{dist}{(x,M_t)}\leq C_{1,0}|s-t|<d$.
Hence $|\xi_s(z)-x|=|\op{dist}{z,M_s}\leq|x-z|<d$ so $z\in C_x$
and $N_t(\xi_s(z))\geq 1$
Since $N_t$ is locally constant we must have $N_t(x)\geq 1$.
Thus $N_t(x)\geq 1$ for all $x\in M_s$.

Let $\bz_s^t=\bz_t\circ\bz_s^{-1}$ for $t\in\cJ(s,\delta_2)$.
\begin{align*}
\int N_t\,d||M_s||
&=\int\cW_m(D\xi_s(z)\circ\tau_t(z))\,d||M_t||z\\
&=\int\cW_m(D\xi_s(\bz_s^t(x)))\cW_m(D\bz_s^t(x))\,d||M_s||x
&\rightarrow \ ||M_s||(\R^n)
\end{align*}
since
$$\lim_{t\rightarrow s}
\sup\{|\cW_m(D\xi_s(\bz_s^t(x)))\cW_m(D\bz_s^t(x))-1|:x\in M_s\}
=0.$$

Thus there is $\delta\in(0,\delta_3]$ such that
$$N_t(x)=1\quad\hbox{if $t\in\cJ(s,\delta)$ and $x\in M_s$.}$$
It follows that $\Lambda_t$ carries $M_s$ diffeomorphically onto $M_t$
for such $t$.

\end{proof}

\subsection{$\ul{M}$.}

We let
$$\ul{M}=\{(s,z):s\in(0,\epsilon_1) \ \hbox{and} \ z\in M_s\}.$$

\begin{thm} $\ul{M}\in\cM_{\bj+m}(\cJ\oplus\R^n)$.
\end{thm}
\begin{proof} Suppose $s\in\cJ(\epsilon_1)$. Let $\delta=\epsilon_1-|s|$. Let
$$\Psi:\cJ(s,\delta)\times N_s(\ubr)\rightarrow\cJ(s,\delta)\times N_s(\ubr)$$
be such that $\Psi(t,z)=(t,\bz_s^{-1}(\xi_s(z)))$ for $(t,z)\in\cJ(s,\delta)\times N_s(\ubr)$ and note that $\Psi$ is smooth.

Let $\delta$ and $\Lambda$ be as in Theorem \ref{Lambda};
note that $\Lambda$ is smooth and let
$$\Gamma:\cJ(s,\delta)\times M\rightarrow\cJ(s,\delta)\times M$$
be such that $\Gamma_t=\Lambda_t^{-1}$ for $t\in\cJ(s,\delta)$ and note that
$\Gamma$ is smooth.

Let
$$\Phi:\cJ(s,\delta)\times M\rightarrow\cJ(s,\delta)\times N_s(\ubr)$$
be such that
$$\Phi(t,p)=(t,\bz_t(\Lambda_t^{-1}(p)))\quad\hbox{for $(t,p)\in\cJ(s,\delta)\times M$.}$$
Note that $\Phi$ is smooth.

If $(t,p)\in\cJ(s,\delta)\times M$ then
\begin{align*}
\Psi(\Phi(t,p))
&=\Psi(t,\bz_t(\Lambda_t^{-1}(p)))\\
&=(t,\bz_t(\bz_s^{-1}(\xi_s(\Lambda_t^{-1}(p)))))\\
&=(t,p).
\end{align*}
The Theorem now follows from Proposition \ref{mani2}.
\end{proof}


\subsection{$\protect\underline{\tau}$, $\protect\underline{\nu}$, $\protect\underline{\gamma}$, $\protect\underline{B}$, $\protect\underline{A}$.}
We define smooth maps
$$\ul{\tau}:\ul{M}\rightarrow\mO_m(\R^n);\quad
\ul{\nu}:\ul{M}\rightarrow\mO_{n-m}(\R^n);\quad
\ul{\gamma}:\ul{M}\rightarrow\mO_{n-m-1}(\R^n)$$
by requiring that
$$\ul{\tau}(s,x)(u)=Di_{\ul{M}}(x,s)(u,0)\quad\hbox{for $(s,x)\in\ul{M}$
and $u\in\R^n$}$$
and letting
$$\ul{\nu}(s,x)=i_{\R^n}-\ul{\tau}(x,s);\quad
\ul{\gamma}(s,x)=\ul{\nu}(s,x)-\iota(x)$$
for $(s,x)\in\ul{M}$ and $u\in\R^n$. It follows that
$$\ul{\tau}(s,x)=\tau_s(x);\quad
\ul{\nu}(s,x)=\nu_s(x);\quad
\ul{\gamma}(s,x)=\gamma_s(x)$$
for $(x,s)\in\ul{M}$.

We define the smooth maps
$$\ul{B}:\ul{M}\rightarrow\MultiLin(\R^n\times\R^n,\R^n)
\quad\hbox{and}\quad
\ul{A}:\ul{M}\rightarrow\Lin(\R^n,\End(\R^n))$$
by letting
$$\ul{B}(s,x)(u,v)=\ul{\nu}(s,x)(D\ul{\tau}(s,x)(0,u)(0,v))
\quad\hbox{whenever $(s,x)\in\ul{M}$ and $u,v\in\R^n$}$$
and requiring that
$$\ul{A}(s,x)(u)(v)\bullet w=u\bullet\ul{B}(s,x)(v,w)
\quad\hbox{whenever $(s,x)\in\ul{M}$ and $u,v,w\in\R^n$.}$$
It follows that
$$\ul{B}(s,x)(u,v)=B_s(x)(u,v)
\quad\hbox{and}\quad
\ul{A}(s,x)(u)=A_s(x)(u)$$
whenever $(s,x)\in\ul{M}$ and $u,v\in\R^n$.

\subsection{$\protect\underline{N}$, $\protect\underline{\Xi}$, $\protect\underline{\xi}$, $\protect\underline{\sigma}$.} 

For each $r\in(0,\ubr]$ we let
$$\ul{N}(r)=\{(s,z)\in \cJ(\epsilon_1)\times\R^n:\op{dist}{(z,M_s)}<r\}.$$
We let
$$\ul{\Xi}:\bE(\ul{\nu},\ubr)\rightarrow\ul{N}(\ubr)$$
be such that
$$\ul{\Xi}((s,x),v)=(s,\Xi_s(x,v))\quad\hbox{for $((s,x),v)\in\bE(\ul{\nu},\ubr)$.}$$
It follows from Theorem \ref{Xidiffeo} that $\ul{\Xi}$ carries $\bE(\ul{\nu},\ubr)$
diffeomorphically onto $\ul{N}(\ubr)$; in particular, $\ul{N}(r)$
is an open subset of $\cJ\times\R^n$ whenever $0<r\leq\ubr$.
Since $\ul{\Xi}$ carries $\bE(\ul{\nu},\ubr)$ diffeomorphically onto $\ul{N}(\ubr)$
we may define smooth functions
$$\ul{\xi}:\ul{N}(\ubr)\rightarrow \R^n\quad\hbox{and}\quad\ul{\sigma}:\ul{N}(\ubr)\rightarrow\R^n$$
by requiring that
$$(\ul{\xi},\ul{\sigma})
=(\ul{\Xi}|\bE(\ul{\nu},\br))^{-1}:\ul{N}(\ubr)\rightarrow\bE(\nu,\br).$$
Note that
\begin{equation}\label{ulxisigma}
(s,z)=\ul{\Xi}(\ul{\xi}(s,z),\ul{\sigma}(s,z))=(s,\ul{\xi}(s,z)+\ul{\sigma}(s,z))
\quad\hbox{for $(s,z)\in \ul{N}(\ubr)$.}
\end{equation}

\subsection{$\protect\underline{f}$, $\protect\underline{\upsilon}$, $\protect\underline{\rho}$, $\protect\underline{V}$, $\protect\underline{g}$.}

We define the smooth maps
$$\ul{f}:\bE(\ul{\nu},\ubr)\rightarrow\cJ\times\R^n;\quad
\ul{\upsilon}:\ul{M}(\ubr)\rightarrow\R^+;\quad
\ul{\rho}:\ul{M}(\ubr)\rightarrow\R^n$$
by letting
$$\ul{f}((t,z),v)=(t,f_t(x,v);\quad
\ul{\upsilon}((t,z),v)=(t,\upsilon_t(x,v);\quad
\ul{\rho}((t,z),v)=(t,\rho_t(x,v)$$
for $((t,x),v))\in\bE(\ul{\nu})$.
We also set
$$\ul{V}=\op{rng}{\ul{f}}.$$

\begin{thm} $\ul{V}$ is an open subset of $\cJ\times\R^n$ and
$\ul{f}$ carries $\bE(\ul{\nu},\ubr)$ diffeomorphically onto $\ul{V}$
with inverse $(\ul{\xi},\ul{\rho})$.

Moreover, $\ul{f}$ carries $\bE(\ul{\gamma},\ubr)$ diffeomorphically
onto $\ul{V}\cap\Sm$.
\end{thm}
\begin{proof} This follows from the results of Section \ref{f}.
\end{proof}

For $r\in(0,\ubr)$ we let
$$\ul{O}(r)=\{(s,t,z)\in\cJ(\epsilon_1)\times\cJ(\epsilon_1)\times\R^n:
(s,z)\in\ul{N}(r) \ \hbox{and} \ (t,z)\in\ul{N}(r)\}.$$
Since $\ul{N}(r)$ is an open subset of $\cJ\times\R^n$ it follows
that
$\ul{O}(r)$ is an open subset of $\cJ\times\cJ\times\R^n$.

\begin{lem}\label{O3} Suppose $s,t\in\cJ(\epsilon_1)$, $q,r\in\R^+$, $q+r\leq\br$
and $(s,t,z)\in\ul{O}(r)$. If $\chi^{1}(\bz)|s-t|\leq q$ then
$$\{(\bc_{s,t}(\lambda),z):\lambda\in[0,1]\}\subset\ul{O}(q+r)$$
\end{lem}
\begin{proof} Since $(s,z)\in\ul{N}(r)$, $\op{dist}{(z,M_s)}<r\leq\ubr\leq\br_s$
so $z\in\op{dmn}{\xi_s}$ which implies that there is $w\in M_s$ such
that $\op{dist}{(z,M_s)}=|z-w|$.
Suppose $\lambda\in[0,1]$ and $u=(1-\lambda)s+\lambda t$.
Then
$$\op{dist}{(z,M_u)}\leq|z-\bz_u(p)|
\leq|z-w|+|\bz_s(p)-\bz_u(p)|
<r+\chi^{1}(\bz)\lambda|t-s|\leq r+q$$
so $z\in\ul{N}(q+r)$. By the same argument with $s$ replaced by $t$
we find that $z\in\ul{N}(q+r)$.
\end{proof} 
\section{ $F_s^t$.} Suppose $s\in \cS$ and $W\in\cE(\gamma_s,\br_s)$.
Let
      $$F_s(W):M_s\rightarrow V_s$$
      be such that
      $$F_s(W)=f_s\circ(i_{M_s},W).$$

Then $\op{rng}{F_s(W)}$ is a smooth compact submanifold of $\Sm$,
$F_s(W)$ carries $M_s$ diffeomorphically onto
$\op{rng}{F_s}(W)$,
$$\xi_s\circ F_s(W)=i_{M_s}
\quad\hbox{and}\quad
\rho_s\circ F_s(W)=0.$$

\begin{prop} Suppose $t\in \cS$, $W_i\in\cE(\gamma_s,\br_s)$, $i=1,2$
and $\op{rng}{F_s(W_1)}=\op{rng}{F_s(W_2)}$. Then $W_1=W_2$.
\end{prop}
\begin{proof} Suppose $x_1\in M_s$. There is $x_2\in M_s$ such that
$F_s(W_2)(x_2)=F_s(W_1)(x_1)$. Thus
$$f_s(x_1,W_1(x_1))=F_s(W_1)(x_1)=F_s(W_2)(x_2)=f_s(x_2,W_2(x_2)).$$
Applying $g_s$ we find that $(x_1,W_1(x_1))=(x_2,W_2(x_2))$.
\end{proof}

\begin{prop} Suppose $s\in \cS$, $X\in\cM_m(\R^n)$ and $X\subset V_s$.
There
is $W\in\cE(\gamma_s,\br_s)$ such that $\op{rng}{F_s(W)}=X$ if and only if
$\xi_s|X$ carries $X$ diffeomorphically onto $M_s$.
\end{prop}
\begin{proof} Suppose $W\in\cE(\gamma_s,\br_s)$ and $\op{rng}{F_s(W)}=X$. Then
$F_s(W)=f_s\circ(i_{M_s},W)$ so
$$(\xi_s\circ F_s(W),\rho_s\circ F_s(W))
=g_s\circ F_s(W)=g_s\circ f_s\circ(i_{M_s},W)=(i_{M_s},W)$$
so $\xi_s\circ F_s(W)=i_{M_s}$. It follows that 
$\xi_s|X$ carries $X$ diffeomorphically onto $M_s$.

Suppose $\xi_s|X$ carries $X$ diffeomorphically onto $M_s$. Let
$W\in\cE(\gamma_s,\br_s)$ be such that
$W\circ(\xi_s|X)=\rho_s|X\in\cE(M_s,\gamma_s)$. Then
$$F_s(W)\circ(\xi_s|X)=f_s\circ(i_{M_s},W)\circ(\xi_s|X)
=f_s\circ(\xi_s|X,\rho_s|X)
=i_X$$
so $\op{rng}{F_s(W)}=X$.
\end{proof}

\begin{defn} For $s,t\in \cS$ we let
$$F_s^t=\{(W,Z)\in\cE(\gamma_s,\br_s)\times\cE(\gamma_t,\br_t):
\op{rng}{F_s(W)}=\op{rng}{F_t(Z)}\}.$$
\end{defn}

\begin{prop} Suppose $s,t\in \cS$, $W\in\cE(\gamma_s,\br_s)$
and $Z\in\cE(\gamma_t,\br_t)$.
The following statements are equivalent.
\begin{itemize}
\item[(I)] $(W,Z)\in F_s^t$;
\item[(II)] $\xi_t\circ F_s(W)$ carries $M_s$ diffeomorphically onto $M_t$
and
\begin{equation}\label{Fsteqn}
Z\circ\xi_t\circ F_s(W)=\rho_t\circ F_s(W).
\end{equation}
\end{itemize}
\end{prop}
\begin{proof} Suppose (I) holds. Let $g=F_t(Z)^{-1}\circ F_s(W)$. Then
$g$ carries $M_s$ diffeomorphically onto $M_t$ and
$F_t(Z)\circ g=F_s(W)$ so for any $y\in M_t$ we have
$$(i_{M_t},Z)\circ g=(\xi_t,\rho_t)\circ F_t(Z)\circ g
=(\xi_t,\rho_t)\circ F_s(W)$$
so $g=\xi_t\circ F_s(W)$ and $Z\circ g=\rho_t\circ F_s(W)$
so (II) holds.


Suppose (II) holds. Let $X=\op{rng}{F_s(W)}$. (\ref{Fsteqn}) implies that
$Z\circ(\xi_t|X)=\rho_t|X$. So
$$F_t(Z)\circ(\xi_t|X)
=f_t\circ(i_{M_t},Z)\circ(\xi_t|X)
=f_t\circ(\xi_t|X,\rho_t|X)
=f_t\circ g_t\circ i_X
=i_X$$
so $\op{rng}{F_t(Z)}=\op{rng}{i_X}=X=\op{rng}{F_s(W)}$.
\end{proof}

\begin{cor}\label{Lst} If $s,t\in\cS$ then $0\in\op{dmn}{F_s^t}$ if and only if
$\xi_t$ carries $M_s$ diffeomorphically onto $M_t$.
\end{cor}

\begin{prop}\label{LipW} Suppose $s\in \cS$ and $W\in\cE(\gamma_s)$. Then
$$\op{Lip}{W}\leq \left(1+{8\over\br_s}\right)\chi^{[1]}(W)
\quad\hbox{and}\quad
\op{Lip}{F_s(W)}\leq 1+\left(1+{8\over\br_s}\right)\chi^{[1]}(W)$$
\end{prop}
\begin{proof} The first inequality follows from
Proposition \ref{Lip}.

Suppose $x_i\in M_s$, $i=1,2$. Using Proposition \ref{Lipbbeta} we find that
\begin{align*}
|F_s(W)(x_1)-F_s(W)(x_2)|
&=|\bbeta(x_1+W(x_2))-\bbeta(x_2+W(x_2))|\\
&\leq|x_1-x_2|+|W(x_1)-W(x_2)|
\end{align*}
so the second inequality holds.
\end{proof}

\begin{prop}\label{Box} Suppose $E$ is a Euclidean space, $h:\ul{N}(\ubr)
\rightarrow E$, $h$ is smooth and
$$\Box h:\ul{O}(\br/2)\rightarrow\Lin(\cJ,E)$$
is such that
$$\Box h(s,t,z)(u)=\int_0^1 D^{1,0}h(\bc_{s,t}(\lambda),z)(u)\,d\lambda
\quad\hbox{for $(s,t,z)\in\ul{O}(\ubr/2)$ and $u\in\cJ$.}$$
Then $\Box h$ is smooth and
\begin{equation}\label{hst}
h(s,z)-h(t,z)=\Box h(s,t,z)(s-t)\quad\hbox{for $(s,t,z)\in\ul{O}(\ubr/2).$}
\end{equation}
If $(s,t,z_i)\in\ul{O}(\ubr/2)$ for $i=1,2$ then
\begin{equation}\label{hst2}
||\Box h(s,t,z_1)-\Box h(s,t,z_2)||\leq C_{1,\Box}(h)|z_1-z_2|
\end{equation}
where we have set
$$C_{1,\Box}(h)=\max\left\{\chi^{2}(h),8{\chi^{0}(h)\over\ubr}\right\}.$$
\end{prop}
\begin{proof} Since $\ul{O}(\br)$ is an open subset
of $\cJ(\epsilon_1)\times\cJ(\epsilon_1)\times\R^n$ it is elementary that $\Box h$ is smooth.
(\ref{hst}) follows since
$$h(s,z)-h(t,z)=\int_0^1{d\over d\lambda}h(\bc_{s,t}(\lambda),z)\,d\lambda
=\int D^{1,0}h(\bc_{s,t}(\lambda),z)(s-t)\,d\lambda$$
whenever $(s,t,z)\in\ul{O}(\ubr/2)$. 

Suppose$(s,t,z_i)\in\ul{O}(\ubr/2)$ for $i=1,2$.
If $|z_2-z_1|<\ubr/4$ then
Corollary \ref{O3} implies that $(s,t,w)\in\ul{O}(\ubr/2)$ if $w\in\ob{z_1}{\ubr/4}$
so that
$$||\Box h(s,t,z_1)-\Box h(s,t,z_2)||\\
\leq\chi^{1}(\Box h)\,|z_1-z_2|\leq\chi^{2}(h)\,|z_1-z_2|.$$
On the other hand, if $|z_1-z_2|\geq\ubr/4$
then
$$||\Box h(s,t,z_1)-\Box h(s,t,z_2)||\leq
2\chi^{1}(h)\leq 8{\chi^{1}(h)\over\ubr}|z_1-z_2|.$$
So (\ref{hst2}) holds.
\end{proof}

\subsection{$\epsilon_2$ and $\cG$.}
Let
$$\epsilon_2\in(0,\epsilon_1]$$
be such that
$$C_{1,\Box}(\ul{\xi})(2\epsilon_2)<{1\over 2}.$$

Let
$$\cG=\left\{(s,W)\in\cJ(\epsilon_1)\times\cE(\gamma_s):
  \left(1+{8\over\ubr}\right)\chi^{[1]}(W)<1 \right\}.
$$
\begin{prop}\label{cGW} Suppose $(s,W)\in\cG$. Then
$$\chi^0(W)<{\ubr\over 8}, \ \op{Lip}{W}<1 \ \hbox{and} \  \op{Lip}{F_s(W)}<2.$$
\end{prop}
\begin{proof} $\chi^0(W)\leq\chi^{[1]}(W)<\ubr/8$. The second and third
inequalities follow directly from Proposition \ref{LipW}.
\end{proof}

\begin{prop} Suppose $s,t\in \cJ(\epsilon_2)$ and $(s,W)\in\cG$.
 Then  
\begin{equation}\label{Ost}
\{(s,t,F_s(W)(x)):x\in M_s\}\subset\ul{O}(\ubr/2);
\end{equation}
in particular, $\op{rng}{F_s(W)}\subset\op{dmn}{\xi_t}$.
\end{prop}
\begin{proof} Suppose $x\in M_s$. There is $p\in M$ such that $\bz_s(p)=x$.
By Proposition \ref{bbetaxv} we have
$$\op{dist}{(F_s(W)(x),M_s)}
\leq|\bbeta(x+W(x))-x|
\leq \chi^0(W)
\leq{\ubr\over 4}$$
and
$$\op{dist}{(F_s(W)(x),M_t)}
\leq|\bbeta(x+W(x))-x|+|\bz(s,p)-\bz(t,p)|
\leq \chi^0(W)+C_{1,0}|s-t|\leq{\ubr\over 2}$$
so (\ref{Ost}) holds and, since $s,t\in\cJ(\epsilon_2)$, $\op{rng}{F_s(W)}
\subset\op{dmn}{\xi_t}$ since $C_{1,\Box}(\ul{\xi})(2\epsilon_2)<1/2$.
\end{proof}

For $(s,W)\in\cG$ and $(t,x)\in \ul{O}(\ubr/2)$ we let
$$A_s^t(W)(x)=\Box \ul{\xi}(s,t,F_s(W)(x))(s-t)$$
and we let
$$B_s^t(W)(x)=\Box \ul{\rho}(s,t,F_s(W)(x))(s-t).$$

\begin{thm} Suppose $s,t\in \cJ(\epsilon_2)$ and $(s,W)\in\cG$.
\begin{equation}\label{xist}
\xi_t(F_s(W)(x))
=x-A_s^t(W)(x);
\end{equation}
\begin{equation}\label{rhost}
\rho_t(F_s(W)(x))
=W(x)-B_s^t(W)(x);
\end{equation}
\begin{equation}\label{biLip}
\begin{split}
{1\over 2}&|x_1-x_2|\\
&\leq(1-2C_{1,\Box}(\ul{\xi})|s-t|)|x_1-x_2|\\
&\leq|\xi_t(F_s(W)(x_1))-\xi_t(F_s(W)(x_2))|\\
&\leq(1+2C_{1,\Box}(\ul{\xi})|s-t|)|x_1-x_2|\\
&\leq{3\over 2}|x_1-x_2|
\end{split}
\end{equation}
if $x_i\in M_s$, $i=1,2$;
and
\begin{equation}\label{inFst}
  W\in\op{dmn}{F_s^t}.
\end{equation}
\end{thm}
\begin{proof}  
(\ref{xist}) and (\ref{rhost}) follow from (\ref{hst})
since $\xi_s(F_s(W)(x))=x$ and $\rho_s(F_s(W)(x))=W(x)$ for $x\in M_s$.

Suppose  $x_i\in M_s$, $i=1,2$. From (\ref{xist}) we infer that
\begin{align*}
\ul{\xi}(t,&F_s(W)(x_1))-\ul{\xi}(t,F_s(W)(x_2))\\
&=x_1-x_2
 -(\Box \ul{\xi}(s,t,F_s(W)(x_1))(s-t)-\Box \ul{\xi}(s,t,F_s(W)(x_2))(s-t))
\end{align*}
and, by Proposition \ref{Box} and Proposition \ref{LipW},
\begin{align*}
|\Box \ul{\xi}(s,t,F_s(W)(x_1))(s-t)&-\Box \ul{\xi}(s,t,F_s(W)(x_2))(s-t)|\\
&\leq C_{1,\Box}(\ul{\xi})|F_s(W)(x_2)-F_s(W)(x_1)||s-t|\\
&\leq C_{1,\Box}(\ul{\xi})2|x_1-x_2||s-t|\\
&\leq {1\over 2}
\end{align*}
It now follows from Proposition \ref{diffeo} that $W\in\op{dmn}{F_s^t}$.
\end{proof}

\begin{thm} Suppose $s,t\in\cJ(\epsilon_2)$ and $(s,W)\in\cG$.
Then
\begin{equation}\label{Dxist}
D(\xi_t\circ F_s(W))(x)
=\tau_s(x)+\Box(D^{0,1}\ul{\xi})(t,s,F_s(W)(x))(t-s)\circ DF_s(W)(x);
\end{equation}
\begin{equation}\label{Drhost}
D(\rho_t\circ F_s(W))(x)
=DW(x)+\Box(D^{0,1}\ul{\rho})(t,s,F_s(W)(x))(t-s)\circ DF_s(W)(x);
\end{equation}
\begin{equation}\label{Fstest}
\chi^0(F_s^t(W))\leq\chi^0(D^{1,0}\ul{\rho})|s-t|+\chi^0(W)
\end{equation}
and
\begin{equation}\label{DFstest}
\chi^0(DF_s^t(W))\leq
\chi^0(D^{1,1}\ul{\rho})|s-t|+\chi^0(DW).
\end{equation}
\end{thm}
\begin{proof} Suppose $x\in M_s$. Using (\ref{hst}) we find that
\begin{align*}
D(\xi_t\circ F_s(W))&(x)-D(\xi_s\circ F_s(W))(x)\\
&=(D^{0,1}\ul{\xi}(t,F_s(W)(x))-D^{0,1}\ul{\xi}(s,F_s(W)(x)))\circ DF_s(W)(x)\\
&=(\Box(D^{0,1}\ul{\xi})(t,s,F_s(W))(t-s))\circ DF_s(W)(x).
\end{align*}
Moreover, $\xi_s\circ F_s(W)=i_{M_s}$ so $D(\xi_s\circ F_s(W))(x)=\tau_s(x)$
so (\ref{Dxist}) holds. With regard to (\ref{Drhost}) we have
$\rho_s\circ F_s(W)=W$ so a similar analysis yields (\ref{Drhost}).

By (\ref{Fsteqn}) we have
$$F_s^t(W)\circ\xi_t\circ F_s(W)=\rho_t\circ F_s(W).$$
(\ref{Fstest}) follows from (\ref{rhost}). The Chain Rule implies that
for any $x\in M_s$ we have
$$DF_s^t(W)((\xi_t\circ F_s(W)(x)))\circ D(\xi_t\circ F_s(W))(x)
=D(\rho_t\circ F_s)(x);$$
(\ref{DFstest}) now follows from (\ref{Dxist}) and (\ref{Drhost}).
\end{proof}





\subsection{$\cL_m$, $\cW_m$ and $\epsilon_3$.}\label{cLcW}
Let $\cL_m=\{L\in\End(\R^n):\op{rank}{L}=m\}$ and let
$$\cW_m:\cL_m\rightarrow\R^+$$
be such that $\cW_m(L)=|\altdown{m}L|$ for $L\in\cL_m$.
It is elementary that if $L\in\cL_m$ and $p$ is orthogonal projection
of $\R^n$ onto the kernel of $L$ then
$$|\cW_m(L)-1|\leq(2^m-1)||L-p||.$$

Suppose $(s,W)\in\cG$ and $(s,t,x)\in\ul{O}(\ubr/2)$.
Let
$$\bc_s^t(W)(x)=\cW_m(D(\xi_t\circ F_s(W))(x))-1.$$
Since $D(\xi_t\circ F_s(W))(x)\in\cL_m$ we infer from (\ref{LipW}) and (\ref{Drhost}) that
\begin{equation}\label{cstest}
\begin{split}
|\bc_s^t(W)(x)|
&\leq(2^m-1)||D(\xi_t\circ F_s(W))(x)-D(\xi_s\circ F_s(W))(x)||\\
&\leq(2^m-1)||D\xi_t(F_s(W)(x))-D\xi_s(F_s(W)(x))||\,||D(F_s(W))(x)||\\
&\leq C_{\cW_m}|s-t|
\end{split}
\end{equation}
where we have set
$$C_{\cW_m}=(2^m-1)2\chi^{1,1}(\ul{\xi}).$$

Let
$$\epsilon_3\in(0,\epsilon_2]$$
be such that
$C_{\cW_m}(2\epsilon_3)<1$. Suppose $s,t\in\cJ(\epsilon_3)$.
Then $-|\bc_s^t(W)|\leq \bc_s^t(W)\leq|\bc_s^t(W)|$ so
\begin{equation}\label{cW1}
 0
 \leq1-C_{\cW_m}|s-t|
 \leq 1-|\bc_s^t(W)|
 \leq 1-\bc_s^t(W)
 \leq 1+|\bc_s^t(W)|
 \leq 1+C_{\cW_m}|s-t|).
\end{equation}

\begin{thm}\label{Fstestimate} Suppose $s,t\in\cJ(\epsilon_3)$
and  $(s,W)\in\cG$. Then
$$\{(s,t,F_s(W)(x)):x\in M_s\}\subset\ul{O}(\ubr/2)$$
and, for any $\epsilon\in\R^+$,
$$|P_t(F_s^t(W),F_s^t(W))-P_s(W,W)|
\leq (\epsilon+C_{F,1}|s-t|)P_s(W,W)+\left(1+{1\over\epsilon}\right)C_{F,2}|s-t|^2$$
where
$$C_{F,1}=C_{\cW_m}\quad\hbox{and}\quad C_{F,2}=\chi^0(\Box\ul{\rho})^2||M_s||(\R^n)(1+C_{\cW_m})(2\epsilon_3).$$
\end{thm}
\begin{proof} Keeping in mind (\ref{cW1}) we have that
\begin{align*}
P_t(F_s^t(W),F_s^t(W))
&=\int|F_s^t(W)|^2\,d||M_t||\\
&=\int |F_s^t(W)\circ\xi_t\circ F_s(W)|^2(1+\bc_s^t(W))\,d||M||_s\\
&=\int |\rho_t\circ F_s(W)|^2(1+\bc_s^t(W))\,d||M||_s\\
&=\int |W-B_s^t(W)|^2(1+\bc_s^t(W))\,d||M||_s\\
&=\int (|W|^2-2W\bullet B_s^t(W)+|B_s^t|^2)(1+\bc_s^t(W))\,d||M||_s
\end{align*}
so, using (\ref{cstest}), we find that 
\begin{align*}|
|P_t(F_s^t(W),&F_s^t(W))-P_s(W,W)|\\
&=\left|\int \bc_s^t(W)|W|^2+(-2W\bullet B_s^t(W)
               +|B_s^t|^2)(1+\bc_s^t(W))\,d||M||_s\right|\\
              &\leq\int(C_{\cW_m}|s-t|+\epsilon)|W|^2\\
              &\quad\quad\quad
                +\left(1+{1\over\epsilon}\right)|B_s^t(W)|^2)(1+C_{\cW_m}|s-t|)\,d||M_s||.
\end{align*}
Since $|B_s^t(W)|^2\leq(\chi^0(\Box\ul{\rho}))^2|s-t|^2$ our estimate holds.
\end{proof}

\section{The family of interest.}

Suppose $M\in\cM_m(\R^n)$, $M\subset\Sm$ and
\begin{equation}\label{H}
H(x)=-x\quad\hbox{for $x\in M$}
\end{equation}
where $H$ is the mean curvature normal of $M$. As is well known and is implied
by Proposition \ref{firstvar}
(\ref{H}) is equivalent to the statement that
\begin{equation}\label{EH}
{d\over d\epsilon}\cH^m(h_t[M])\big|_{\epsilon=0}=0
\end{equation}
whenever $I$ is an open interval $I$ containing $0$, $h:I\times M\rightarrow\Sm$,
$h$ is smooth and $h_0(x)=x$ for $x\in M$.

Let $\gamma$ be as in \ref{gamma} and let $\br$ be the reach of $M$.

\subsection{The symmetric quadratic functionals $P,Q,R,S$.}
We define symmetric quadratic functionals
$$ P,Q,R,S:\cE(\gamma)\times\cE(\gamma)\rightarrow\R$$
by letting
\begin{align*}
P(X,Y)&=\int X\bullet Y\,d||M||;\\
Q(X,Y)&=\int D_\gamma X\bullet D_\gamma Y\,d||M||;\\
R(X,Y)&=\int A(x)(X(x))\bullet A(x)(Y(x))\,d||M||x;\\
S(X,Y)&=Q(X,Y)-m P(X,Y)-R(X,Y)
\end{align*}
for $X,Y\in\cE(\gamma)$.

\subsection{$\Phi$; $\Psi$; $\bL$; and $E_\lambda$, $\lambda\in\R$.}
For $Z\in\cE(\gamma,\br)$ we let $F(Z)$ be as in (\ref{F})
and we let
$$\Phi(Z)=\cH^m(\op{rng}{F(Z)}).$$

We define
$$\Psi:\cC^{2,\alpha}(\gamma,\br)\rightarrow\cC^{0,\alpha}(\gamma)$$
by requiring that
$$P(\Psi(Z),W)={d\over d\epsilon}\Phi(Z+\epsilon W)\big|_{\epsilon=0}.$$
(In fact, the right hand side of this formula is an integral over $M$ of linear function of $DW$  and the preceding formula
follows by an integration by parts in the integral to obtain an integral over
$M$ of a linear function of $W$.)
Let
$$\bL\in\Lin(\cE(\gamma),\cE(\gamma))$$
be the linearization at $0$ of $\Psi$; that is,
$$\bL(Z)={d\over d\eta}\Psi(\eta Z)\big|_{\eta=0}
\quad\hbox{for $Z\in\cE(\gamma)$.}$$
As is well known,
\begin{equation}\label{symmetric}
P(\bL(Z),W)=P(Z,\bL(W))=S(Z,W)\quad\hbox{whenever $Z,W\in\cE(\gamma)$.}
\end{equation}

For each $\lambda\in\R$ we let
$$E_\lambda=\{Z\in\cE(\gamma):\bL(Z)=\lambda Z\};$$
so $Z\in E_{\lambda}$ if and only if
$$S(Z,W)=\lambda P(Z,W)\quad\hbox{for $W\in\cE(\gamma)$.}$$

\subsection{The H\"older spaces $\cC^{k,\alpha}(\gamma)$.}
Let us fix $\alpha\in(0,1)$.
For each $k\in\N$ we let $\cC^{k,\alpha}(\gamma)$ be the vector space
of maps $Z:M\rightarrow\R^n$ such that $\gamma(x)(Z(x))=Z(x)$ for each
$x\in M$; $D^j Z$ exists and is continuous for each $j\in\N$ such that $j\leq k$; and $D^k Z$ is H\"older continuous with exponent $\alpha$. $\cC^{k,\alpha}(\gamma)$ has a natural Banach space structure.
Clearly, $P$ and $R$ have uniqued continuous extensions, still denoted by $P$
and $R$, to $\cC^{0,\alpha}(\gamma)$
and $Q$ and $S$ have natural continuous extensions, still denoted by $Q$
and $S$, to $\cC^{1,\alpha}(\gamma)$. $\bL$ has an unique extension,
still denoted by $\bL$ to a continuous map from $\cC^{2,\alpha}(\gamma)$
to $\cC^{0,\alpha}(\gamma)$. Moreover, (\ref{symmetric}) holds when
$Z,W\in\cC^{2,\alpha}$.

$\Phi$ has a unique continuous extension to $\{Z\in\cC^{1,\alpha}(\gamma):
|Z|<\br\}$ and $\Psi$ has a unique extension to a continuous map
from $\{Z\in\cC^{2,\alpha}(\gamma)\}$ to $\cC^{0,\alpha}(\gamma)$.
These extensions are real analytic.

Let
$$\cJ=\op{ker}{\bL} \ \hbox{and let} \ \bj=\op{dim}{\cJ}.$$
The members of $\cJ$ are called {\bf Jacobi fields relative to $\Sm$}.
{\em The inner product on $\cJ$ is the restriction to $\cJ\times\cJ$
of $P$.} So if $s,t\in\cJ$ then
$$s\bullet t=\int s(x)\bullet t(x)\,d||M||.$$
By standard elliptic theory,
$$\bj \ \hbox{is finite;}$$
for each $\epsilon\in\R^+$ there is a constant $C$ such that
$$\sum_{|\lambda|\leq N}\op{dim}{E_\lambda}\leq(1+N)^{m/2+\epsilon}
\quad\hbox{whenever $N\in\R^+$;}$$
in particular, $\op{dim}{\op{rng}{E_\lambda}}$ is finite for each
$\lambda\in\R$ and $\{\lambda\in\R:E_\lambda\neq 0\}$ is discrete;
for each $s\in\R$ and each $Z\in\cE(\gamma)$,
$$\Sigma_s(Z)=P(E_0(Z),E_0(Z))+\sum_{\lambda\in\R\sim\{0\}}|\lambda|^2 P(E_\lambda(Z),E_\lambda(Z))<\infty$$
and $\{\Sigma_s:s\in\R\}$ is a family of seminorms which define the
strong topology on $\cE(\gamma)$.

\subsection{Construction of the family of interest.}
For $k=0,2$ let $\cH_k$ be the $P$ orthogonal complement
of $\cJ$ in $\cC^{k,\alpha}(\gamma)$; since $\bj$ is finite $\cH_k$
is closed; let $\Pi_k$ be the 
$P$ orthogonal projection of $\cC^{k,\alpha}(\gamma)$
onto $\cJ$ and let $\Pi_k^\perp$ be the complementary projection; thus
$Z=\Pi_k(Z)+\Pi_k^\perp(Z)$ for $Z\in\cC^{k,\alpha}(\gamma)$.
By standard elliptic theory,
$$\hbox{$\bL|\cH_2$ is Banach space isomorphism from $\cH_2$ onto $\cH_0$.}$$
Since $\Psi(0)=0$ and since $\Psi$ is real analytic
it follows from the Implicit Function Theorem that there are open neighborhoods $\cS$ of $0$ in $\cJ$ and of $0$ in $\cV$ of $0$ in $\op{rng}{\Pi_0}^\perp\tilde{\Gamma}$ and a
real analytic function $\Upsilon:\cS\rightarrow \op{rng}{\Pi_0}^\perp$ such that
 $\{s+\Upsilon(s):s\in\cS\}\subset \cV$
and such that if $s\in\cS$, $W\in\cJ^\perp$ and $s+W\in \cV$ then
$$\Pi_0^\perp(\Psi(s+W))=0 \ \Leftrightarrow  \ W=\Upsilon(s).$$
We let $\Omega:\cS\rightarrow\cJ$ be such that
$$\Omega(s)=\Pi_0(\Psi(s+\Upsilon(s))\quad\hbox{for $s\in\cS$}$$
and we let
$$\cA=\{s\in\cS:\Omega(s)=0\}.$$
Note that $\Omega$ is a real analytic function so $\cA$ is a real analytic set.
Evidently,
$$\{Z\in\cV:\Psi(Z)=0\}=\{s+\Upsilon(s):s\in\cS \ \hbox{and} \ \Omega(s)=0\}.$$
Since $\Psi(0)=0$ we have
$$\Upsilon(0)=0\quad\hbox{and}\quad\Omega(0)=0.$$

{\it We now assume that}
\begin{equation}\label{main}
  \Tan{\cA,0}=\cJ;
\end{equation}
owing to the real analyticity of $\cA$ this implies that
$$\Omega \ \hbox{vanishes near $0$ in $\cS$.}$$
{\em That $\Tan{\cA,0}=\cJ$ is clearly implied by the assumption about Jacobi
fields made in (1) of the Introduction of \cite{AA81}.} So we have
$$\{Z\in\cV:\Psi(Z)=0\}=\{s+\Upsilon(s):s\in\cS\}.$$

We define
$$ \bz:\cS\times M\rightarrow\Sm$$
by letting
$$\bz(s,x)=F(s+\Upsilon(s))(x)=f(x,(s+\Upsilon(s))(x))\in\Sm
\quad\hbox{for $(s,x)\in\cS\times M$.}$$
Letting
$$M_s=\op{rng}{F(s+\Upsilon(s))}$$
we find that
$$M_s \ \hbox{is $\Phi$ stationary.}$$
It follows from Proposition \ref{firstvar} that
\begin{equation}\label{stationary}
\cH^m(M_s)=\cH^m(M)\quad\hbox{for $s$ in the connected component of $0$ in $\cS$.}
\end{equation}
\begin{rem} In fact (\ref{stationary}) holds {\em without} assuming (\ref{main})
owing to the fact that the real analyticity of $\cA$ implies that $\cA$
is locally rectifiably arcwise connected.
\end{rem}
Let
$$\cJ_s$$
be the Jacobi fields on $M_s$ relative to $\Sm$.

For $s\in\cS$ we let
$$\gamma_s:M_s\rightarrow\mO_m(\R^n)$$
be such that $\gamma_s(y)$ equals orthogonal projection of $\R^n$ onto
$\Nor{M_s,y}\cap\Tan{\Sm,y}$.
For each $s\in\cS$ let
$$\bJ_s\in\Lin(\cJ,\cE(\gamma_s))$$
be such that
$$\bJ_s(t)(\bz(s,x))
=\gamma_s(\bz(s,x))\left({d\over d\epsilon}\bz(s+\epsilon t,x)\big|_{\epsilon=0}\right)\quad\hbox{for $t\in\cJ$.}$$
Since $M_{s+\epsilon t}$ is stationary for sufficiently small $\epsilon$
we find that
$$\op{rng}{\bJ_s}\subset\cJ_s.$$

\begin{thm} $D\Upsilon(0)=0$ and $\bJ_0(t)=t$ for $t\in\cJ$.
Moreover,
\begin{equation}\label{kerbJ}
s\in\cS \ \Rightarrow \ \op{ker}{\bJ_s}=\{0\}.
\end{equation}
\end{thm}
\begin{proof}
Suppose $s\in\cS$, $t\in\cJ$ and
$W=s+\Upsilon(s)\in\cE(\gamma_s)$. Let $X=t+D\Upsilon(s)(t)\in\cE(\gamma_s)$.
Suppose $x\in M$ and $y=\bz(s,x)$. Using Proposition \ref{Dbbeta}) we
find that
\begin{align*}
\bJ_s&(t)(y)\\
&=\gamma_s(y)\left(
{d\over d\epsilon}(F(s+\epsilon t))(x)\big|_{\epsilon=0}\right)\\
&=\gamma_s(y)\left({d\over d\epsilon}\bbeta(x+(s+\epsilon t+\Upsilon(s+\epsilon t))(x))\big|_{\epsilon=0}\right)\\
&=\gamma_s(y)D\bbeta(x+W(x))(X(x))\\
 &=\gamma_s(y)\left(1+|W(x)|)^{-3/2}\left((1+|W(x)|^2)X(x)-(X(x)\bullet W(x))(x+W(x))\right)\right)\\
  &=(1+|W(x)|)^{-3/2}\left((1+|W(x)|^2)X(x)-(X(x)\bullet W(x))(W(x))\right).\\
  \end{align*}

In particular, if $s=0$ then $W=0$ so that
$$\bJ_0(t)=X=t+D\Upsilon(0)(t).$$
Since $\bJ_0(t)\in\cJ$ and $D\Upsilon(0)(t)\in\cH_0$
we find that $D\Upsilon(0)(t)=0$ and $\bJ_0(t)=t$.

Now suppose $t\in\op{ker}{\bJ_s}$. Then
$$X(x)={(X(x)\bullet W(x))W(x)\over 1+|W(x)|^2}$$
which implies that
$$|X(x)|\leq|X(x)|{|W(x)|^2\over 1+|W(x)|^2}$$
which in turn implies that $X(x)=0$ which implies that $t(x)=0$.
\end{proof}

\begin{cor}\label{bJ0} Suppose $s,t\in\cJ$. Then
$$s\bullet t=P_0(\bJ_0(s),\bJ_0(t)).$$
\end{cor}
\begin{proof} We have
$$P_0(\bJ_0(s),\bJ_0(t))
=\int\bJ_0(s)(x)\bullet\bJ_0(t)(x)\,d||M_0||x
=\int s(x)\bullet t(x)\,d||M||=s\bullet t.$$
\end{proof}

\begin{thm} Suppose $t\in\cJ(\epsilon_1)$, $u\in\cJ$ and $x\in M_t$. Then
$$\bJ_t(u)(x)=-D\ul{\rho}(t,x)(u,0).$$
\end{thm}
\begin{rem} This is 5.4(7) of \cite{AA81}.
\end{rem}
\begin{proof} Let $x=\bz(t,p)\in M_t$.
Let $I=\{\epsilon\in\R:t+\epsilon u\in\cJ\}$  and note that $I$ is an open neigborhood of $0$ in $\cS$.
For $\epsilon\in I$ let
$$g(\epsilon)=\ul{f}(t+\epsilon u,(\bz(t+\epsilon u,p),0))
=\bbeta(\bz(t+\epsilon u,p)).$$
We have $g(0)=\bz(t,p)=x$ and, by Proposition \ref{Dbbeta},
$$g^\prime(0)=\bo(\bz(t,p))(D\bz(t,p)(u,0))=D\bz(t,p)(u,0)$$
since $|\bz(t,p)|=1$ for $(t,p)\in\cS\times M$.
Since $g(\epsilon)\in M_{t+\epsilon u}$ we have
$$0=\ul{\rho}(t+\epsilon u,g(\epsilon))\quad\hbox{for $\epsilon\in I$;}$$
Differentiating at $\epsilon=0$ we find using Proposition \ref{Drho} that
\begin{align*}
0&={d\over d\epsilon}\ul{\rho}(t+\epsilon u,g(\epsilon))\big|_{\epsilon=0}\\
&=D\ul{\rho}(t,\bz(t,p))(u,g^\prime(0))\\
&=D\ul{\rho}(t,\bz(t,p))(u,D\bz(t,p)(u,0))\\
&=D\ul{\rho}(t,\bz(t,p))(u,0)+\gamma_t(\bz(t,p))(D\bz(t,p)(u,0)).
\end{align*}
\end{proof}

\begin{lem} There is a constant $C_{\bJ,0}$ such that
$$|\bJ_s(u)(x)-\bJ_t(u)(\xi_t(x))|\leq C\,|s-t|\,|u|$$
Suppose $s,t\in\cJ(\epsilon_1)$, $u\in\cJ$ and $x\in M_s$. Then
\end{lem}
\begin{proof}
\begin{align*}
\bJ_s(u)&(x)-\bJ_t(u)(\xi_t(x))\\
&=(D\ul{\rho}(s,x)-D\ul{\rho}(t,\xi_t(x)))(u,0)\\
&=(D\ul{\rho}(s,\xi_s(x))-D\ul{\rho}(t,\xi_t(x)))(u,0)\\
&=(D\ul{\rho}(s,\xi_s(x))-D\ul{\rho}(s,\xi_t(x)))(u,0)
+(D\ul{\rho}(s,\xi_t(x))-D\ul{\rho}(t,\xi_t(x)))(u,0).
\end{align*}
Now apply Proposition \ref{Box}.
\end{proof}
\begin{cor}\label{Jst} There is a constant $C_\bJ$ such that
$$|P_s(\bJ_s(u),\bJ_s(v))-P_t(\bJ_t(u),\bJ_t(v))|
\leq C_\bJ|s-t|\,|u|\,|v|\quad\hbox{for $s,t\in\cJ(\epsilon_1)$
and $u,v\in\cJ$.}$$
\end{cor}
\begin{proof} Suppose $s,t\in\cJ(\epsilon_1)$ and $u,v\in\cJ$. For each
$x\in M_s$ let
\begin{align*}
\Gamma_1(x)&=\bJ_s(u)(x);\\
\Gamma_2(x)&=\bJ_s(v)(x);\\
\Gamma_3(x)&=\bJ_t(u)(\xi_t(x));\\
\Gamma_4(x)&=\bJ_t(v)(\xi_t(x));\\
\Gamma_5(x)&=\bJ_t(u)(\xi_t(x))\bullet\bJ_t(v)(\xi_t(x))\bc_s^t(0)(x)).
\end{align*}
Then
\begin{align*}
|P_s(\bJ_s(u)&,\bJ_s(v))-P_t(\bJ_t(u),\bJ_t(v))|\\
&=\int\Gamma_1\bullet\Gamma_2-\Gamma_3\bullet\Gamma_4
-\Gamma_5\,d||M_s||\\
&=\int(\Gamma_1-\Gamma_2)\bullet\Gamma_3
+\Gamma_3\bullet(\Gamma_2-\Gamma_4)-\Gamma_5\,d||M_s||.
\end{align*}
Now apply the preceding Lemmma and (\ref{cstest}).
\end{proof}

The following Corollary follows by letting $t=0$ above.
\begin{cor}\label{Jst2} If $u\in\cJ$ and $s\in\cJ(\epsilon_1)$ then
  $$P_s(\bJ_s(u),\bJ_s(u))\leq (1+C_\bJ)|s|\,|u|.$$
\end{cor}

\subsection{Symmetric bilinear functions.}
Suppose $V$ is an inner product space, $\bB$ is the vector space of
real value symmetric bilinear functions on $V\times V$. We say $T\in\bB$
is {\bf nonnegative} if $T(v,v)\geq 0$ for $v\in V$.and $\bB^+$
is the set of those $T\in\bB$ such that $T(v,v)\geq 0$ for $v\in V$.
Then
\begin{equation}\label{expand}
T(v,w)={1\over 2}(T(v+w,v+w)-T(v,v)-T(w,w)
\end{equation}
and, if $T$ is nonnegative,
\begin{equation}\label{CaS}
|T(v,w)|\leq\sqrt{T(v,v)T(w,w)}\leq{1\over 2}(T(v,v)+T(w,w))
\end{equation}
for $v,w\in V$.

\begin{prop}\label{Tineq} Suppose $T_i\in\bB$, $i=1,2$, $U\in\bB$, $U$ is nonnegative
and
$$|T_1(v,v)-T_2(v,v)|\leq U(v,v)\quad\hbox{for $v\in V$.}$$
Then
$$|T_1(v,w)-T_2(v,w)|\leq{3\over 2}(U(u,u)+U(v,v))\quad\hbox{for $u,v\in V$.}$$
\end{prop}
\begin{proof} Suppose $v,w\in V$. Then

\begin{align*}
|T_1(v,w)&-T_2(v,w)|\\
&={1\over 2}|(T_1(v+w,v+w)\\
&\quad\quad-T_1(v,v)-T_1(w,w))-(T_2(v+w,v+w)-T_2(v,v)-T_2(w,w))|\\
&={1\over 2}|(T_1(v+w,v+w)\\
&\quad\quad-T_2(v+w,v+w))-(T_1(v)-T_2(v))-(T_1(w)-T_2(w))|\\
&\leq{1\over 2}(U(v+w,v+w)+U(v,v)+U(w,w))\\
&=U(v,v)+U(w,w)+U(v,w)\\
&\leq{3\over 2}(U(v,v)+U(w,w)),
\end{align*}
\end{proof}

\section{$L_s^t$ and $\epsilon_4$.}

Suppose $s,t\in\cJ(\epsilon_2)$. We define a linear map
$$L_s^t:\cE(\gamma_s)\rightarrow\cE(\gamma_t)$$
as follows. Keeping in mind that Proposition \ref{Lst} implies that
$\xi_t$ carries $M_s$ diffeomorphically onto $M_t$ we require that
\begin{equation}\label{Lst1}
L_s^t(W)(\xi_t(x))=\gamma_t(\xi_t(x))(W(x))\quad\hbox{for
  $W\in\cE(\gamma_s)$ and $x\in M_s$.}
\end{equation}
\begin{rem} This is {\em not} the $L_s^t$ of Part Two of the proof of
Lemma 5.4 in \cite{AA81}. There
we see
$$L_s^t(W)={d\over d\epsilon}F_s^t(\epsilon W)\big|_{\epsilon=0}.$$
This has the pleasant property that $L_s^u=L_t^u\circ L_s^t$ which $L_s^t$
above does not have but it is harder to work with.
\end{rem}
We let
$$h_L:\ul{N}(\ubr)\rightarrow\End(\R^n)$$
be such that $h_L(s,z)=\gamma_s(\xi_s(z))$ 
for $(s,z)\in \ul{N}(\ubr)$; note that $h_L$ is smooth. We infer from
Proposition \ref{Box} that
$$\gamma_s(\xi_s(z))-\gamma_t(\xi_t(z))
=\Box h_L(s,t,z)(s-t)
\quad\hbox{for $(s,t,z)\in\ul{O}(\ubr/2)$.}$$
If $W\in\cE(\gamma_s)$ and $x\in M_s$ we have
\begin{equation}\label{Lst2}
L_s^t(W)(\xi_t(x))=W(x)-\Box h_L(s,t,\xi_t(x))(s-t)(W(x))
\end{equation}
whenever $x\in M_s$ since
$\gamma_s(\xi_s(x))(W(x))=\gamma_s(x)(W(x))=W(x)$.
Let
$$C_L=\sup\{||\Box h_L(s,t,z)(u)||:(s,t,z)\in \ul{O}(\ubr), \ u\in\cJ
\ \hbox{and} \ |u|\leq 1\}.$$
Then $C_L<\infty$ and
\begin{equation}\label{CL0}
||\gamma_s(\xi_s(z))-\gamma_t(\xi_t(z))||\leq C_L|s-t|
\quad\hbox{whenever $z\in\ul{O}(\ubr/2)$}
\end{equation}
and
\begin{equation}\label{CL}
(1-C_L|s-t|)|W(x)|\leq|L_s^t(W)(\xi_t(x))|\leq(1+C_L|s-t|)|W(x)|.
\end{equation}

Let
$$\epsilon_4\in(0,\min\{1,\epsilon_3\}]$$
be such that $2C_L(2\epsilon_4)<1$.

The following Theorem is an immediate consequence of (\ref{CL}).
\begin{thm}\label{Liso}  Suppose $s,t\in\cJ(\epsilon_4)$.
Then $L_s^t$ carries $\cE(\gamma_s)$ isomorphically onto $\cE(\gamma_t)$Su
\end{thm}

\begin{thm}\label{LP}
and there is a constant $C_{L,P}$ such that
$$|P_t(L_s^t(X),L_s^t(Y))-P_s(X,Y)|\leq C_{L,P}|s-t|)P_s(X,Y).$$ 
\end{thm}
\begin{proof} Suppose $W\in\cE(\gamma_s)$ and $\bc_s^t$ is as in (\ref{cLcW}).Using
\ref{Lst2} we find that
\begin{align*}\int|L_s^t(W)|^2&-|W|^2\,d||M_t||\\
&=\int |L_s^t(W)\circ\xi_t|^2\cW_m(D(\xi_t|M_s)) -|W|^2\,d||M_s||\\
&=\int |L_s^t(W)\circ\xi_t|^2 (1+\bc_s^t(0))\,d||M_s||\\
&=\int |W(x)-\Box h_L(s,t,\xi_t(x)(s-t)|^2(1+\bc_s^t(0)(x))-|W(x)|^2\,d||M_s||\\
\end{align*}
Now use (\ref{cW1}) and Proposition \ref{Tineq} to complete the proof.
\end{proof}

\begin{prop}\label{RleqP} Suppose $s\in\cJ(\epsilon_1)$ and $W\in\cE(\gamma_s)$. Then
$$|R_s(W,W)|\leq m\ubr^2 P_s(W,W).$$
\end{prop}
\begin{proof} Suppose $u_1,\ldots,u_n$ is an orthonormal basis for $\R^n$.
For any $x\in M_s$ we have
$$|A_s(x)(W(x))|^2
=\sum_{i=1}^m |A_s(x)(W(x))(u_i)|^2
\leq m\,||B_s(x)||^2\,||W(x)|^2
\leq{m\ubb}^2|W(x)|^2.$$
\end{proof}

\begin{thm}\label{LR} Suppose $s,t\in\cJ(\epsilon_4)$.
There is a constant $C_{L,R}$ such that
$$|R_t(L_s^t(X),L_s^t(Y))-R_s(X,Y)|\leq C_{L,R}|s-t|)P_s(X,Y)
\quad\hbox{for $X,Y\in\cE(\gamma_s)$.}$$
\end{thm}
\begin{proof} We let
$$h_A(s,z)=\ul{A}(s,\ul{\xi}(s,z))(\ul{\gamma}(s,\ul{\xi}(s,z)))
\quad\hbox{for $(s,z)\in\ul{O}(\ubr)$.}$$

Suppose $W\in\cE(\gamma_s)$.
Since $(s,t,x)\in\ul{O}(\ubr/2)$ we have
$$A_s(x)(W(x))=A_s(\xi_s(x))(\gamma_s(\xi_s(x))(W(x))=h_A(s,x)(W(x))$$
as well as
\begin{align*}
A_t(\xi_t(x))&(L_s^t(W)(\xi_t(x)))\\
&=A_t(\xi_t(x))(\gamma_t(\xi_t(x))(W(x)))\\
&=\ul{A}(t,\ul{\xi}(t,x))(\ul{\gamma}(t,\ul{\xi}(t,x))(W(x)))\\
&=h_A(t,x)(W(x))
\end{align*}
so
\begin{align*}
R_t(L_s^t&(W),L_s^t(W))-R_s(W,W)\\
&=\int|A_t(y)(L_s^t(W)(y))|^2\,d||M_t||y-\int|A_s(x)(W(x)\,d||M_s||x\\
&=\int|A_t(\xi_t(x))(L_s^t(W)(\xi_t(x))|^2(1+\bc_s^t(0)(x))\,d||M_s||x\\
&=\int h_A(t,x)(W(x))(1+\bc_s^t(0)(x))-h_A(s,x)(W(x))\,d||M_s||\\
&=\int (h_A(t,x)-h_A(s,x))(W(x))+h_A(t,x)(W(x))\bc_s^t(0)(x))\,d||M_s||x.
\end{align*}
The existence of the constant $C_{L,R}$ now follows from Proposition \ref{Box},
(\ref{cstest}) and Proposition \ref{Tineq}.
\end{proof}

\begin{thm}\label{LQ} Suppose $s,t\in\cJ(\epsilon_4)$. There is a constant $C_{L,Q}$
such that for $X,Y\in\cE(\gamma_s,\ubr)$ we have
\begin{align*}
|Q_t(L_s^t(X),L_s^t(Y))-Q_s(X,Y)|
\leq C_{L,Q}|s-t|\big((P_t+Q_t)(X,X)+(P_t+Q_t)(Y,Y)\big)
\end{align*}
\end{thm}

\begin{proof} For $u\in\cJ(\epsilon_4)$ and $U,V\in\cE(\gamma_u,\ubr)$ we set
$$T_u(U,V)=\int DU\bullet DV\,d||M_u||$$
and observe that
\begin{align*}
T_u(U,V)
&=\int (D_\gamma U(y)+A_u(y)(U(y)))\bullet(D_\gamma V(y)+A_u(y)(V(y)))\,d||M_u||y\\
&=\int D_\gamma U(y)\bullet D_\gamma V(y)+A_u(y)(U(y))\bullet A_u(y)(V(y))\,d||M_u||y\\
&=Q_u(U,V)+R_u(U,V).
\end{align*}

Suppose $W\in\cE(\gamma_s,\ubr)$. We calculate
\begin{align*}
Q_t(L_s^t(W)&,L_s^t(W))\\
&=\int DL_s^t(W)(y)\bullet DL_s^t(W)(y)\,d||M_t||y\\
&=\int DL_s^t(W)(\xi_t(x))\bullet DL_s^t(W)(\xi_t(x))
       (1+\bc_s^t(x))\,d||M_s||x.
\end{align*}

Now by (ii) of Theorem \ref{Dxisigma} and (\ref{Lst2}) we have
for any $W\in\cE(\gamma_s,\ubr)$ that
\begin{equation}\label{QT}
\begin{split}
DL_s^t(W)&(\xi_t(x))(u)\\
&=DL_s^t(W)(\xi_t(x))\circ\tau_t(\xi_t(x))(u)\\
&=DL_s^t(W)(\xi_t(x))\circ D\xi_t(\xi_t(x))(u)\\
&=D(L_s^t(W)\circ\xi_t)(\xi_t(x))(u)\\
&=DW(x)(u)-D^{0,0,1}\Box h_L(s,t,\xi_t(x))(u)(s-t)(W(x))\\
&\quad\quad\quad-\Box h_L(s,t,\xi_t(x))(s-t)(DW(x)(u))\\
&=\sum_{i=0}^2\Gamma_i(s,t,x)(DW(x),s-t,W(x))(u)
\end{split}
\end{equation}
where
$$\Gamma_i\in\cE(\ul{O}(\ubr/2),
	\MultiLin(\End(\R^n)\times\cJ\times\R^n,\End(\R^n)))
	\quad i=0,1,2,$$
are such that
\begin{align*}
\Gamma_0(s,t,x)(l,v,w)(u)&=l(u);\\
\Gamma_1(s,t,x)(l,v,w)(u)&=-D^{0,0,1}\Box h_L(s,t,\xi_t(x))(u)(v)(w);\\
\Gamma_2(s,t,x)(l,v,w)(u)&=-\Box h_L(s,t,\xi_t(x))(v)(l(u))
\end{align*}

For $i,j=0,1,2$ we let
$$Z_{i,j}(W,W)\in\cE(\ul{O}(\ubr/2)$$ be such that
$$Z_{i,j}(W,W)(s,t,x)
=\Gamma_i(s,t,x)(DW(x),s-t,W(x))\bullet\Gamma_j(s,t,x)(DW(x),s-t,W(x))$$
for $i,j=0,1,2$ and $(s,t,x)\in\ul{O}(\ubr/2)$. Then (\ref{QT}) implies that
\begin{align*}
T_t(L_s^t(W),L_s^t(W))
&=\int DL_s^t(W)(\xi_t(x))\bullet DL_s^t(W)(\xi_t(x))(1+\bc_s^t(0))\,d||M_s||x\\
&=\int \sum_{i,j=0}^2 Z_{i,j}(W,W)(s,t,x)(1+\bc_s^t(0)(x))\,d||M_s||x.
\end{align*}

We have
$$|Z_{0,0}(W,W)-|DW|^2|=|\bc_s^t||DW|^2\leq C_{\cW_m}|s-t||DW|^2.$$
If $i,j\in\{0,1,2\}$ and not both $i$ and $j$ are $0$ there is a constant
$c_1$ such that
$$|Z_{i,j}(W,W)(s,t,x)(1+\bc_s^t(x))|\leq c_1|s-t|(|W|^2+|DW|^2).$$
It follows that
$$|T_t(L_s^t(W),L_s^t(W))-T_s(W,W)|\leq
(C_{\cW_m}+5c_1)C_{\cW_m}|s-t|(P_s(W,W)+T_s(W,W)).$$
Since $Q_s=T_s-R_s$ the Theorem now follows from Theorem \ref{LR}
and Proposition \ref{Tineq}.
\end{proof}

\begin{thm}\label{LS} Suppose $s,t\in\cJ(\epsilon_4)$. There is a constant
$C_{L,S}$ such that
$$|S_t(L_s^t(X),L_s^t(Y))-S_s(X,Y)|\leq C_{L,S}|s-t|((P_s+Q_s)(X,X)+(P_s+Q_s)(Y,Y))\quad\hbox{for $X,Y\in\cE(\gamma_s)$.}$$
\end{thm}
\begin{proof} Suppose $X,Y\in\cE(\gamma_s)$. Then
\begin{align*}
|S_t(L_s^t(X)&,L_s^t(Y))-S_s(X,Y)|\\
&=|(Q_t-R_t-mP_t)(L_s^t(X),L_s^t(Y))
-(Q_s-R_s-mP_s)(X,Y)|\\
&\leq|Q_t(L_s^t(X),L_s^t(Y))-Q_s(X,Y)|\\
&+|R_t(L_s^t(X),L_s^t(Y))-R_s(X,Y)|\\
&+m|P_t(L_s^t(X),L_s^t(Y))-P_s(X,Y)|
\end{align*}
Now apply Theorems \ref{LP}, \ref{LR} and \ref{LQ}.
\end{proof}

l b

\section{$\epsilon_5$.}

Let
$$\epsilon_5\in(0,\epsilon_4]$$
be such that
\begin{equation}\label{epsilon5}
2C_{L,Q}\epsilon_5<{1\over 2}.
\end{equation}

\begin{rem} The following Lemma is a replacement for Part Four of Lemma 5.4 of \cite{AA81}
which is false; to see this let $\lambda$ there equal zero. Here is a corrected version for which the inequalities derived in the proof of Part Four of
Lemma 5.4 of \cite{AA81} do not suffice to prove the replacement.
\end{rem}

\begin{lem}\label{PartFour} There is a constant $C_{L,*}$ such that if $s,t\in\cJ(\epsilon_5)$,
 $X\in\cE(\gamma_s)$, $X^*=X-E_{s,0}(X)$ and
$$S_t(L_s^t(X),Y)=\lambda P_t(L_s^t(X),Y)\quad\hbox{for $Y\in\cE(\gamma_t)$}$$
then
$$\kappa_s P_s(X^*,X^*)
\leq\big(|\lambda|(1+C_{L,*}|s-t|)+C_{L,*}|s-t|\big)P_s(X,X).$$
\end{lem}
\begin{proof} We have
$$\kappa_s P_s(X^*,X^*)\leq |S_s(X,X^*)|\leq\Gamma_1+\Gamma_2$$
where
$$\Gamma_1=|S_t(L_s^t(X),L_s^t(X^*))|
\quad\hbox{and}\quad
\Gamma_2=|S_t(L_s^t(X),L_s^t(X^*))-S_s(X,X^*)|.$$

\begin{align*}
\Gamma_1&=|\lambda||P_t(L_s^t(X),L_s^t(X^*))|\\
&\leq|\lambda|(|P_s(X,X^*)+C_{L,P}|s-t|(P_s(X,X)+P_s(X^*,X^*)))\\
&\leq|\lambda|(1+2C_{L,P}|s-t|)P_s(X,X)
\end{align*}
and
\begin{align*}
\Gamma_2
&\leq C_{L,S}|s-t|\left((P_s+Q_s)(X,X)+(P_s+Q_s)(X^*,X^*)\right)\\
&\leq C_{L,S}|s-t|(2P_s(X,X)+Q_s(X,X)+Q_s(X^*,X^*)).
\end{align*}

We estimate
\begin{align*}
Q_s(&X,X)\\
&\leq Q_t(L_s^t(X),L_s^t(X))+C_{L,Q}|s-t|(P_s(X,X)+Q_s(X,X))\\
&=(S_t+mP_t+R_t)(L_s^t(X),L_s^t(X)+C_{L,Q}|s-t|(P_s(X,X)+Q_s(X,X))\\
&=\lambda P_t(L_s^t(X),L_s^t(X)+(mP_t+R_t)(L_s^t(X),L_s^t(X))\\
&\quad\quad+C_{L,Q}|s-t|(P_s(X,X)+Q_s(X,X))\\
&\leq((|\lambda|+m)(1+C_{L,P}|s-t|)+(1+C_{L,R}|s-t|))P_s(X,X)\\
&\quad\quad+ C_{L,Q}|s-t|(P_s(X,X)+Q_s(X,X))
\end{align*}
so, as $1-C_{L,Q}|s-t|\geq 1-C_{L,Q}(2\epsilon_5)\geq 1/2$ we have
$$Q_s(X,X)\leq 2((|\lambda|+m)(1+C_{L,P}|s-t|)+(1+C_{L,R}|s-t|))P_s(X,X).$$

Since $S_s(X_0,Y)=0$ for any $Y\in\cE(\gamma_s)$ we estimate
\begin{align*}
Q_s(X^*,X^*)&=Q_s(X,X)-2Q_s(X,X_0)+Q_s(X_0,X_0)\\
&=Q_s(X,X)+mP_s(X,X_0)+R_s(X,X_0)+mP_s(X_0,X_0)+R_s(X_0,X_0)\\
&\leq Q_s(X,X)+(2m+2m^2\ubb^2)P_s(X,X).
\end{align*}
The existence of $C_{L,*}$ follows.
\end{proof}

\subsection{$\Pi_s$.}

For each $s\in\cJ(\epsilon_5)$ let
$$\Pi_s$$
be $Pi_s$ orthogonal projection of $\cE(\gamma_s)$ onto $\op{rng}{\bJ_s}$.

\begin{lem} There is a constant $C_\Pi$ such that
$$|P_t(\Pi_t(L_s^t(W),\Pi_t(L_s^t(W)))-P_s(\Pi_s(W),\Pi_s(W))|=
C_\Pi|s-t|P_s(W,W)$$
whenever $s\in\cJ(\epsilon_5)$ and $W\in\cE(\gamma_s)$.
\end{lem}
\begin{proof}
Let $A_s=\bJ_s^*\circ\bJ_s$. Then $\bJ_s\circ A_s^{-1/2}$ is an isometry
from $\cJ$ to $\op{rng}{\bJ_s}$. Indeed, if $u,v\in\cJ$ we have
\begin{equation}\label{Aiso}
\begin{split}
P_s(\bJ_s\circ A_s^{-1/2}(u),&\bullet\bJ_s\circ A_s^{-1/2}(v))\\
&=\bJ_s^*\circ\bJ_s\circ A_s^{-1/2}(u)\bullet v\\
&=A_s\circ A_s^{-1/2}(u)\bullet A_s^{-1/2}(v)\\
&=u\bullet v.
\end{split}
\end{equation}

Suppose $W\in\cE(\gamma_s)$.
If $u_1,\ldots,u_{\bj}$ is an orthonormal basis for $\cJ$ and
$U_{s,i}=\bJ_s^*\circ A^{-1/2}(u_i)\in\cE(\gamma_s)$ then
$$\Pi_s(W)=\sum_{i=1}^\bj P_s(W,U_{s,i})U_{s,i}$$ 
which implies that
$$\Pi_t(L_s^t(W),L_s^t(W))=\sum_{i=1}^\bj P_t(L_s^t(W),U_{t,i})^2.$$
Thus
\begin{align*}
P_t(L_s^t(W),U_{t,i})
&=\int L_s^t(W)\bullet U_{t,i}\,d||M_t||\\
&=\int L_s^t(W)(\xi_t(x))\bullet U_{t,i}(\xi_t(x))(1+\bc_s^t(x))\,d||M_s||x\\
&=\int \gamma_t(\xi_t(x))(W(x))\bullet U_{t,i}(\xi_t(X))(1+\bc_s^t(x))\,d||M_s||x
\end{align*}
The Lemma now follows from the estimates occurring in the proof of Theorem \ref{LR}.
\end{proof}

\subsection{$\epsilon_6$.}
\begin{rem} The following Lemma is, essentially, Part Six of Lemma 5.4 of
\cite{AA81}.
\end{rem}

\begin{lem}\label{PartSix} Suppose $s,t\in\cJ(\epsilon_5)$, $\Pi_s=E_{s,0}$, $\lambda\in\R$ 
and $\Pi_t\neq E_{t,\lambda}$. Then
$$\kappa_s(1-C_{L,P}|s-t|)\leq |\lambda|(1+C_{L,*}|s-t|)+C_{L,*}|s-t|.$$
\end{lem}
\begin{proof} Choose $Y$ such that $Y\in\op{rng}{E_{t,\lambda}}\sim\{0\}$ and
$\Pi_t(Y)=0$; if $\lambda=0$ this is possible since $\op{rng}{\Pi_t}\subset
E_{t,0}$ and if $\lambda\neq 0$ this is possible since $\op{rng}{\Pi_t}\subset
E_{t,0}$ which does not meet $E_{t,\lambda}$. Let $X\in\cE(\gamma_s)\sim\{0\}$ be
such that $Y=L_s^t(X)$; this is possible since $L_s^t$ carries $\cE(\gamma_s)$
isomorphically onto $\cE(\gamma_t)$.
Rescaling $Y$ if necessary may assume that $P_s(X,X)=1$.

We have
\begin{align*}
P_s(\Pi_s(X),\Pi_s(X))
&=|P_t(Y,Y)-P_s(\Pi_s(X),\Pi_s(X))|\\
&=|P_t(L_s^t(X),L_s^t(X))-P_s(\Pi_s(X),\Pi_s(X))|\\
&\leq C_{L,P}|s-t|P_s(X,X)\\
&=C_{L,P}|s-t|
\end{align*}
 so that, as $\op{rng}{\Pi_s}\subset \op{rng}{E_{s,0}}$,
\begin{align*}
\kappa_s(1&-C_{L,P}|s-t|)
&\leq\kappa_s(1-P_s(\Pi_s(X),\Pi_s(X)))\\
&\leq\kappa_s(P_s(X,X)-P_s(E_{s,0}(X),E_{s,0}(X)));
\end{align*}
now apply Lemma \ref{PartFour}.
\end{proof}

Let
$$\epsilon_6\in(0,\epsilon_5]$$
be such that $\kappa_0>(C_{L,P}+C_{L,*})\epsilon_6$.

\begin{cor} $\op{rng}{\bJ_t}=\op{rng}{E_{t,0}}$ for $t\in\cJ(\epsilon_6)$.
\end{cor}
\begin{proof} Suppose $t\in\cJ(\epsilon_6)$ and $\Pi_t\neq E_{t,0}$.
Applying Lemma \ref{PartSix} with $\lambda=0$ we would have
$\kappa_0(1-C_{L,P}|t|)\leq C_{L,*}|t|$ which would imply  $\kappa_0\leq(C_{L,S}(1+C)|t|<(C_{L,P}+C_{L,*})\epsilon_6$.
\end{proof}

\begin{cor} There is a constant $C_\kappa$ such that
$$s,t\in\cJ(\epsilon_6) \ \Rightarrow \ |\kappa_s-\kappa_t|\leq C_\kappa|s-t|.$$
\end{cor}
\begin{proof} For $u\in\cJ(\epsilon_6)$ let $\Lambda_u=\{\lambda\in\R\sim\{0\}:E_{t,\lambda}\neq 0\}$.

Suppose $u\in\cJ(\epsilon_6)$ and $\lambda\in\Lambda_u$. Since $\Pi_0=E_{0,0}$ and $\Pi_u=E_{u,0}\neq E_{u,\lambda}$  Lemma \ref{PartSix} implies
$$\kappa_u(1-C_{L,P}|u|)\leq|\lambda|(1+C_{L,*}|u|)+C_{L,*}|u|.$$
Since $\kappa_0=\inf\Lambda_0$ we have
$$\kappa_u(1-C_{L,P}|u|)\leq\kappa_0(1+C_{L,*}|u|)+C_{L,*}|u|.$$

Now suppose $s,t\in\cJ(\epsilon_6)$, $\kappa_t\geq\kappa_s$ and $\lambda\in\Lambda_s$. Since $\Pi_s=E_{s,\lambda}$, Lemma \ref{PartSix} implies that
$$\kappa_t(1-C_{L,P}|s-t|)\leq|\lambda|(1+C_{L,*}|s-t|)+C_{L,*}|s-t|.$$
Since $\kappa_s=\inf\Lambda_s$ we find that
$$\kappa_t(1-C_{L,P}|s-t|)\leq\kappa_s(1+C_{L,*}|s-t|)+C_{L,*}|s-t|;$$
subtracting $\kappa_s(1-C_{L,P}|s-t|)$ from both sides we obtain
\begin{align*}
  (\kappa_t-\kappa_s)&(1-C_{L,P}|s-t|)\\
  &\leq\kappa_s(C_{L,*}+C_{L,P})|s-t|+C_{L,*}|s-t|\\
                                 &\leq{\kappa_0(1+C_{L,*}|s|)+C_{L,*}|s|\over 1-C_{L,P}|s|}
                                   (C_{L,*}+C_{L,P})|s-t|+C_{L,*}|s-t|
\end{align*}
so the Corollary holds.

\end{proof}

\section{Proof of a statement similar to 5.4(4) of  \cite{AA81}.}

\begin{thm} There is $\delta\in\R^+$ such that if $(s,Z)\in\cG$ and $|s|<\epsilon_6/2$ and $\chi^{[1]}(Z)<\delta$ there is $t\in\cJ(s,\epsilon_6/2)$ such
that
$$|s-t|^2\leq 2(1+C_\bJ|s|)P_s(E_{s,0}(Z),E_{s,0}(Z))$$
and
$$ E_{t,0}(F_s^t(Z))=0.$$
\end{thm}

\subsection{$\zeta_{Z,s}$.}

For $(s,Z)\in\cG$ with $|s|<\epsilon_6/2$ we define
$$\zeta_{Z,s}:\cJ(\epsilon_6)\rightarrow\cJ$$
by requiring that
$$\zeta_{Z,s}(t)\bullet v
=\int F_s^t(Z)\bullet \bJ_t(v)\,d||M_t||
\quad\hbox{for $t\in \cJ(\epsilon_6)$ and $v\in\cJ$.}$$

Since $F_s^s(Z)=Z$ we infer from Corollary \ref{Jst2} that
\begin{align*}
|\zeta_{Z,s}(s)\bullet v|
&=\leq\chi^0(Z)\int|\bJ_s(v)|\,d||M_s||\\
&\leq\chi^0(Z)P_s(\bJ_s(v),\bJ_s(v))^{1/2}\\
&\leq\chi^0(Z)(1+C_\bJ|s|)^{1/2}|v|
\end{align*}
so
\begin{equation}\label{gs}
|\zeta_{Z,s}(s)|\leq{\epsilon_6\over 4} \quad\hbox{if}\quad
\chi^0(Z)\left(1+C_\bJ{\epsilon_6\over 2}\right)^{1/2}<{\epsilon_6\over 4}.
\end{equation}

\begin{lem}\label{DzetaZ0} There is a constant $C_{1,\zeta}$ such that
$$||D\zeta_{Z,s}(t)-D\zeta_{0,s}(t)||\leq C_{1,\zeta}\chi^{[1]}(Z).$$
\end{lem}
\begin{proof} We have
\begin{align*}
\zeta_{Z,s}(t)\bullet v
&=\int F_s^t(Z)(\xi_t(F_s(Z)(x)))\bullet \bJ_t(v)(\xi_t(F_s(Z)(x)))
\cW_m(D(\xi_t\circ F_s(Z))(x))\,d||M_s||\\
 &=\int \rho_t(F_s(Z)(x)\bullet \bJ_t(v)(\xi_t(F_s(Z)(x)))
    \cW_m(D(\xi_t\circ F_s(Z)\,d||M_s||x\\
&=\int A_s(Z)(t,x)B_s(Z)(t,x)(v)C_s(Z)(t,x)\,d||M_s||x
\end{align*}
where
\begin{align*}
A_s(Z)(t,x)&=\rho_t(F_s(Z)(x))=\ul{\rho}(t,F_s(Z)(x));\\ 
B_s(Z)(t,x)(v)&=\bJ_t(v)(\xi_t(F_s(Z)(x)))=-D^{1,0}\ul{\rho}(t,\ul{\xi}(t,F_s(Z)(x)))(v);\\ 
C_s(Z)(t,x)&=\cW_m(D(\xi_t\circ F_s(Z))(x))
=\cW_m(D^{0,1}\ul{\xi}(t,F_s(Z)(x))\circ DF_s(Z)(x));
\end{align*}

Consequently,
\begin{align*}
D^{1,0}A(Z)(t,x)(u)
&=D^{1,0}\ul{\rho}(t,F_s(Z)(x))(u);\\
D^{1,0}B(Z)(t,x)(u)(v)
&=- D(D^{1,0}\ul{\rho})(t,\ul{\xi}(t,F_s(Z)(x)))(u,D^{1,0}\xi(t,F_s(Z)(x))(u))(v);\\
D^{1,0}C(W)(t,x)(u)
&=D\cW_m(D^{0,1}\ul{\xi}(t,F_s(Z)(x))\circ DF_s(Z)(x)))\\
&\quad\quad\quad(D^{1,0}(D^{0,1}\ul{\xi})(t,F_s(Z)(x)(u)\circ DF_s(Z)(x)))).
\end{align*}
In particular,
\begin{align*}
D^{1,0}A(0)(t,x)(u)
&=D^{1,0}\ul{\rho}(t,x)(u);\\
D^{1,0}B(0)(t,x)(u)(v)
&=- D(D^{1,0}\ul{\rho})(t,\ul{\xi}(t,x))(u,D^{1,0}\xi(t,x)(u))(v);\\
D^{1,0}C(W)(t,x)(u)
&=D\cW_m(D^{0,1}\ul{\xi}(t,x)\circ \tau_s(x))\\
&\quad\quad\quad(D^{1,0}(D^{0,1}(\ul{\xi})(t,x)(u)\circ \tau_x(x))).
\end{align*}

\begin{lem} There are constants $c_A,c_B,c_C$ such that
$$||D^{1,0}A(Z)(t,x)-D^{1,0}A(0)(t,x)||\leq c_A\chi^{[1]}(Z);$$
$$||D^{1,0}B(Z)(t,x)-D^{1,0}B(0)(t,x)||\leq c_B\chi^{[1]}(Z);$$
$$||D^{1,0}C(Z)(t,x)-D^{1,0}C(0)(t,x)||\leq c_C\chi^{[1]}(Z).$$
\end{lem}
\begin{proof} Observe that $F_s(Z)(x)=f_s(x,Z(x))$ and $F_s(0)(x)=f_s(x,0)$.
This implies that $DF_s(Z)(x)=Df_s(x,Z(x))\circ(\tau_s(x),DZ(x))$
and $DF_s(0)(x)=Df_s(x,0)\circ(\tau_s(x),0)$.
\end{proof}

\end{proof}

\begin{lem}\label{zeta2} There is a constant $C_{2,\zeta}$ such that
$$||D\zeta_{0,s}(t)-D\zeta_{0,s}(s)||\leq C_{2,\zeta}|s-t|
\quad\hbox{for $s\in\cJ(\epsilon_6)$ and $t\in\cJ(\epsilon_6)$.}$$
\end{lem}
\begin{proof}
We have
\begin{align*}
\zeta_{0,s}(t)\bullet v
&=\int \rho_t(x)\bullet D^{1,0}\ul{\rho}(t,x)(v)
\cW_m(D\xi_t(x)\circ\tau_s(x))\,d||M_s||.\\
\end{align*}
It follows that
\begin{align*}
D\zeta_{0,s}&(t)(u)\bullet v\\
&=\int D^{1,0}\ul{\rho}(t,x)(u)\bullet D^{1,0}\ul{\rho}(t,x)(v)
      \cW_m(D\xi_t(x)\circ\tau_s(x))\\
&\quad +\rho_t(x)\bullet D^{2,0}\ul{\rho}(t,x)(u\odot v)\cW_m(D\xi_t(x)\circ \tau_s(x))\\
&\quad +\rho_t(x)\bullet D^{1,0}\ul{\rho}(t,x)(v)D\cW_m(D\xi_t(x)\circ\tau_s(x))(D^{1,0}\ul{\xi}(t,x)(u)\circ\tau(x)\,d||M_s||x
\end{align*}

Since $\rho_s(x)=0$ and $D\xi_s(x)\circ\tau_s(x)=\tau_s(x)$ for $x\in M_s$ we have
\begin{equation}
D\zeta_{0,s}(s)(u)\bullet v
=\int D^{1,0}\ul{\rho}(s,x)(u)\bullet D^{1,0}\ul{\rho}(s,x)(v)\,d||M_s||x
\end{equation}
If $(t,x)\in\ul{N}(\ubr)$ then $\ul{\rho}(t,x)=-\Box\ul{\rho}(s,t,x)(s-t)$
and
$$D^{1,0}\ul{\rho}(t,x)
=D^{1,0}\ul{\rho}(s,x)-\Box (D^{1,0}\ul{\rho})(s,t,x)(s-t).$$
\end{proof}

\begin{lem}\label{zeta3} We have
$$||D\zeta_{0,s}(s)-i_{\cJ}||\leq C_\bJ|s|
\quad\hbox{for $s\in\cJ(\epsilon_1)$.}$$
\end{lem}
\begin{proof} Suppose $s\in\cJ(\epsilon_6)$ and $u,v\in\cJ$. Then
\begin{align*}
D\zeta_{0,s}&(s)(u)\bullet v-u\bullet v\\
&=\int D^{1,0}\ul{\rho}(s,x)(u)\bullet D^{1,0}\ul{\rho}(s,x)(v)\,d||M_s||x
       -\int\bJ_0(u)\bullet\bJ_0(v)\,d||M_0||\\
&=\int \bJ_s(u)\bullet \bJ_s(v)\,d||M_s||
-\bJ_0(u)\bullet\bJ_0(v)\,d||M_0||\\
&=P_s(\bJ_s(u),\bJ_s(v))-P_0(\bJ_0(u),\bJ_0(v));
\end{align*}
Now apply Corollary \ref{Jst}.
\end{proof}

\subsection{Finding a zero of $\zeta_{Z,s}$.}
Let $B=\{u\in\cJ:|u-s|\leq\zeta_6/2\}\subset\cJ(\epsilon_6)$; let $g:B\rightarrow\cJ$
be such that $g(u)=\zeta_{Z,s}(u)$ for $u\in B$; and let

It follows from Lemmas (\ref{DzetaZ0}), (\ref{zeta2}) \ref{zeta3} that
$$||D\zeta_{Z,s}(t)-i_{\cJ}||
\leq||D\zeta_{Z,s}(t)-D\zeta_{0,s}(s)||+||D\zeta_{0,s}(s)-i_{\cJ}||
\leq C_{1,\zeta}\chi^{[1]}(Z)+C_{2,\zeta}|s-t|+C_\bJ|s|.$$
Suppose $(\ref{gs})$; then
$$|g(s)|<{\epsilon_6\over 4}.$$
Suppose
\begin{equation}\label{Lipp}
C_{1,\zeta}\chi^{[1]}(Z)+C_{2,\zeta}|s-t|+C_\bJ|s|\leq{1\over 2}.
\end{equation}
$p:B\rightarrow\cJ$ be such that $p(u)=(g(u)-g(s))-(u-s)$
for $u\in B$. Then $p(s)=0$ and $\op{Lip}{p}\leq 1/2$.
Let $C:B\rightarrow\cJ$ be such that
$C(u)=s-g(s)-p(u)$ for $u\in B$. Then $\op{Lip}{C}=\op{Lip}{p}\leq 1/2$.

Suppose $u\in B$. Then
$$|C(u)-s|=|g(s)+p(u)|=|g(s)+p(u)-p(s)|
\leq |g(s)|+\op{Lip}{p}|u-s|\leq\epsilon_6-|s|$$
so $C[B]\subset B$. Since $B$ is closed and $\op{Lip}{C}\leq 1/2<1$, $C$ has a unique fixed point $t\in B$. Since
$t=C(t)=s-g(s)-p(t)=s-g(s)-(p(t)-p(s)$ we have
$$|s-t|=|g(s)+(p(t)-p(s))|\leq |g(s)|+{|s-t|\over 2}$$
so
$$|s-t|\leq 2|g(s)|.$$

Suppose $u\in\cJ$ and $|u|=1$. Then
$$g(s)\bullet u
=\int Z\bullet \bJ_s(u)\,d||M_s||
=P_s(Z,\bJ_s(u))
=P_s(E_{s,0}(Z),\bJ_s(u))$$
so by Corollary \ref{Jst} we have
\begin{align*}
|g(s)\bullet u|^2&\leq P_s(E_{s,0}(Z),E_{s,0}(Z))P_s(\bJ_s(u),\bJ_s(u))\\
&\leq P_s(E_{s,0}(Z),E_{s,0}(Z))(P_0(\bJ_0(u),\bJ_0(u))+C_\bJ|s|)\\
&\leq P_s(E_{s,0}(Z),E_{s,0}(Z))(1+C_\bJ|s|)
\end{align*}
so that
$$|s-t|^2\leq2|g(s)|\leq P_s(E_{s,0}(Z),E_{s,0}(Z))(1+C_\bJ|s|).$$

\end{document}